\documentclass[11pt]{amsart}
\usepackage{amsmath, amsthm, amssymb, amscd, amsfonts}
\usepackage[all]{xy}
\usepackage{newlfont}
\usepackage[dvips]{graphics}
\usepackage[latin1]{inputenc}
\usepackage{eucal}
\usepackage{latexsym}
\usepackage[dvips, dvipsnames, usenames]{color}
\input xy
\xyoption{all}
\usepackage{ifthen}

\newcommand{\VGamma}{\widehat{\Gamma}}

\def\de{\Delta}

\def\bdem{\begin{proof}}
\def\edem{\end{proof}}
\def\N{\mathbb{N}}

\def\G{\mathbb{G}}

\def\cC{\mathcal{C}}

\def\cG{\mathcal{G}}
\def\cH{\mathcal{H}}
\def\cK{\mathcal{K}}

\def\cS{\mathcal{S}}

\def\cW{\mathcal{W}}

\def\cX{\mathcal{X}}

\def\Z{\mathbb{Z}}

\def\bB{\mathfrak{B}}

\def\bS{\mathfrak{S}}

\def\bq{\mathfrak{q}}
\def\cG{\mathcal{G}}
\def\cR{\mathcal{R}}

\def\kk{\mathsf{k}}

\def\yd{{}^{H}_{H}\mathcal{YD}}

\def\xx{\mathbb{X}}

\def\unon{\left\{ 1,\ldots ,\theta \right\}}
\def\zt{\Z^{\theta}}

\newcommand\id{{\operatorname{id}}}
\newcommand\ad{\operatorname{ad}}

\newcommand\ord{{\operatorname{ord}}}
\newcommand\Hom{{\operatorname{Hom}}}
\newcommand\Aut{{\operatorname{Aut}}}
\newcommand\Sh{{\operatorname{Sh}}}

\newcommand{\ot}{{\otimes}}
\newcommand{\otv}{{1\le i \le \theta}}
\newcommand{\otvz}{{1\le i, j \le \theta}}

\newlength{\mpb}

\numberwithin{equation}{section} \theoremstyle{plain}
\newtheorem{theorem}{Theorem}[section]
\newtheorem{cor}[theorem]{Corollary}
\newtheorem{lema}[theorem]{Lemma}

\newtheorem{prop}[theorem]{Proposition}

\theoremstyle{definition}
\newtheorem{defn}[theorem]{Definition}

\newtheorem{question}[theorem]{Question}
\newtheorem{conjecture}[theorem]{Conjecture}
\newtheorem{obs}[theorem]{Remark}

\def\bp{\begin{proof}}
\def\ep{\end{proof}}

\begin{document}

\title[A presentation of Nichols algebras of diagonal type]{A presentation by generators and relations of Nichols algebras of diagonal type and convex orders on root systems}

\author{Iv\'an Ezequiel Angiono}
\address{Facultad of Matem\'atica, Astronom\'i a y F\'isica
\newline \indent
Universidad Nacional of C\'ordoba
\newline
\indent CIEM -- CONICET
\newline
\indent (5000) Ciudad Universitaria, C\'ordoba, Argentina}
\email{angiono@famaf.unc.edu.ar}

\date{\today}


\begin{abstract}
We obtain a presentation by generators and relations of any Nichols algebra of diagonal type with finite root system. We prove that the defining ideal is finitely generated. The proof is based in Kharchenko's theory of PBW basis of Lyndon words. We prove that the lexicographic order on Lyndon words is convex for such PBW generators and so the PBW basis is orthogonal with respect to the canonical non-degenerate form associated to the Nichols algebra.
\end{abstract}

\setlength{\unitlength}{1mm} \settowidth{\mpb}{$q_0\in k^\ast
\setminus \{-1,1\}$,}

\maketitle

\section*{Introduction}

The consideration of pointed Hopf algebras has grown since the appearance of quantized enveloping algebras \cite{Dr,Ji}.
The finite-dimensional analogues, the so-called small quantum groups, were introduced and described by Lusztig \cite{L2,L3}.

The Lifting Method of Andruskiewitsch and Schneider is the leading method for the classification of finite-dimensional
pointed Hopf algebras. Such method depends on the answers to some questions, including the following one:

\smallskip

\begin{question}
\cite[Question 5.9]{A}: Given a braided vector space of diagonal type, determine if the associated Nichols algebra is
finite-dimensional, and in such case compute its dimension. Give a nice presentation by generators and relations.
\end{question}

The first part of this question has been answered by Heckenberger in \cite{H2}, where the author gives a list of all diagonal braidings
whose associated Nichols algebra has a finite root system, but neither an explicit formula for the dimension nor a finite set of defining relations are given. Some of them are Lusztig's examples, which are associated with
the so-called Cartan braidings and for which the dimension and a presentation by generators and relations are known.
Standard braidings were introduced in \cite{AA} and they constitute a family which includes properly the family of
Cartan braidings. Nichols algebras with standard braidings have been presented by generators and relations in \cite{An}, where also an explicit formula for the dimension has been given. Another result about presentation of Nichols algebras is given in \cite{Y} for quantized enveloping algebras associated with semisimple Lie superalgebras. Some other preliminaries considerations on the relations of a Nichols algebra of diagonal type appear in \cite{He}.

Andruskiewitsch and Schneider \cite{AS2} have classified finite-dimensional pointed Hopf algebras
whose group of group-like elements is abelian of order not divisible by some small primes using the Lifting method; all the possible such braidings are of finite Cartan type. They
answered positively the following conjecture for $H_0=\kk \Gamma$, $\Gamma$ an abelian group as above:

\begin{conjecture} \cite[Conj. 1.4]{AS} \label{Conj:Nico} Let $H$ be a finite-dimensional pointed Hopf algebra over $\kk$. Then $H$ is generated by group-like and skew-primitive elements.
\end{conjecture}

This result was proved as a previous step of the main Theorem in \cite{AS2} using the presentation by generators and relations. The conjecture was recently proved in
a more general context \cite{AnGa}, when the braiding is of standard type. The proof follows also using
the presentation by generators and relations.

Because of the braidings of Cartan type we see that there exists a close relation between pointed Hopf algebras and the classical Lie theory.
In such direction the definition of the Weyl groupoid and the root system \cite{H1, HS, HY} associated to a Nichols algebra
$\bB(V)$ of diagonal type has shown to be a good extension of the idea of root systems and Weyl groups associated to semisimple
Lie algebras. Such root system is obtained as the set of degrees of the generators of any PBW basis, and controls coideal subalgebras
between other structures associated to $\bB(V)$ \cite{HS}.

In the classical case, convex orders over the root system were described in order to characterize quantized enveloping algebras $U_q(\mathfrak{g})$ for $\mathfrak{g}$ semisimple \cite{KT,Le, R2}, and to obtain Lusztig isomorphisms in the affine case \cite{Be}. The characterization of convex orders is in consequence necessary, and it has been done for finite \cite{P} and affine \cite{I} root systems.

Our main result is Theorem \ref{Thm:presentacion}: we obtain a presentation by generators and relations for any Nichols algebra of diagonal type whose root system is finite. We obtain two kind of relations that are enough to present $\bB(V)$: powers of root vectors (generators of a PBW basis), and some generalizations of quantum Serre relations which express the braided bracket of two root vectors as a linear combination of other root vectors in an explicit way, see Lemmata \ref{Lemma:heigthgenerators}, \ref{Lemma:otrasrelaciones}.

Theorem \ref{Thm:presentacion} follows by consideration of PBW bases as in \cite{Kh}. Such PBW bases consist of homogeneous polynomials
associated to Lyndon letters (which are called hyperletters) and inherit the lexicographical order. Another important element is the characterization of convex orders for generalized root systems. Such convex orders are related with reduced expressions of elements of the Weyl groupoid. These reduced expressions characterize also right coideal subalgebras of Nichols algebras, so we can relate convex orders and coideal subalgebras.
In particular, the following result holds by Theorem \ref{Thm:presentacion}:

\begin{theorem}
 Let $V$ be a braided vector space of diagonal type whose associated root system in finite, and let $I(V)$ be the ideal of $T(V)$
such that $\bB(V)= T(V)/I(V)$. Then $I(V)$ is finitely generated.
\end{theorem}

Theorem \ref{Thm:presentacion} extends the presentation by generators and relations of Nichols algebras of standard type, see Remark \ref{Rem:standard}, and then gives a new proof for braidings of Cartan type. In particular we obtain the classical presentation of quantized enveloping algebras $U_q(\mathfrak{g})$ and Lusztig's small quantum groups $\mathfrak{u}_q(\mathfrak{g})$, with a different proof. We hope that this presentation helps to prove Conjecture \ref{Conj:Nico} when the group of group-like elements is abelian, see Remark \ref{Rem:generacion en rango uno}.

The plan of this article is the following. In Section \ref{section:preliminaries} we recall the definition of Nichols algebra. We also consider results from \cite{Kh, R2} concerning a PBW basis for Nichols algebras of diagonal type.

In Section \ref{section:rootsystems-coidealsubalgebras} we deal with root systems and coideal subalgebras of Nichols algebras
of diagonal type. In Subsection \ref{subsection:weylgroupoid} we recall the notion of Weyl groupoid and root system, and give
some properties of these objects. In Subsection \ref{subsetion:convexorder} we characterize convex orders on finite root systems
generalizing the results in \cite{P}. In Subsection \ref{subsection:coidealsubalgebras} we recall some results from \cite{HS}
involving coideal subalgebras of Nichols algebras of diagonal type with finite root systems and use these results to characterize
PBW bases of hyperletters. In particular we obtain that the lexicograpical order on the hyperletters is convex.

In Section \ref{section:relations} we obtain the desired presentation by generators and relations. First we prove that the
Kharchenko's PBW basis is orthogonal for the canonical non-degenerate bilinear form as in Proposition  \ref{Prop:formabilineal}
when the braiding matrix is symmetric. Power root vector relations hold in $\bB(V)$ by Lemma \ref{Lemma:heigthgenerators} and generalized quantum Serre relations hold by Lemma \ref{Lemma:otrasrelaciones}. These two sets of relations are enough to give the presentation. We show
in Section \ref{section:examples} how the main theorem allows to obtain explicitly the presentation of some Nichols algebras in some examples.

\medbreak \textbf{Notation.} $\N$ denotes the set of positive integers, and $\N_0$ the set of non-negative integers.

We fix an algebraically closed field $\kk$ of characteristic 0; all vector spaces, Hopf algebras and
tensor products are considered over $\kk$.

For each $N > 0$, $\G_N$ denotes the group of $N$-th roots
of 1 in $\kk$.

Given $n \in \N$, we set the following polynomials in $q$:
$$ \binom{n}{j}_q = \frac{(n)_q!}{(k)_q! (n-k)_q!}, \quad \mbox{where }(n)_q!= \prod_{j=1}^n (k)_q, \quad \mbox{and } (k)_q= \sum_{j=0}^{k-1} q^j. $$

\section{Preliminaries}\label{section:preliminaries}

We recall some definitions and results that we shall need in the subsequent sections. They are related with characterizations of Nichols algebras of diagonal type and PBW bases of such algebras.

Recall that a braided vector space is a pair $(V,c)$, where $V$ is a vector
space and $c\in \Aut (V\ot V)$ is a solution of the braid
equation: $$(c\otimes \id) (\id\otimes c) (c\otimes \id) =
(\id\otimes c) (c\otimes \id) (\id\otimes c).$$
A braided vector space $(V,c)$ is of \emph{diagonal type} if there exists a basis $x_1, \dots x_\theta$ and scalars
$q_{ij}\in \kk^{\times}$ such that
\begin{equation}\label{ecuacion de trenzas}
c(x_i \ot x_j)= q _{ij} x_j \ot x_i, \quad \otvz.
\end{equation}
Following \cite{Kh} we describe an appropriate PBW-basis of a braided graded Hopf algebra $\bB = \oplus_{n\in \N} \bB^n$ such that $\bB^1 \cong V$, where $(V,c)$ is of diagonal type. In particular we obtain PBW bases for Nichols algebras $\bB(V)$ of diagonal type. This construction is based in the notion of Lyndon words. Each Lyndon word has a canonical decomposition as a product of a
pair of smaller Lyndon words, called the Shirshov decomposition. Using such decomposition and the braided bracket, we define inductively a set of hyperwords, which are the elements of a PBW basis for braided graded Hopf algebras of diagonal type. We recall also some properties of this PBW basis.

\subsection{Braided vector spaces of diagonal type and Nichols algebras}\label{subsection:bvs}

Given a braided vector space $(V,c)$, this braiding can be extended to $c:T(V)\ot T(V) \to T(V)\ot T(V)$ canonically, see \eqref{braiding} for the diagonal case. We define for each pair
$x,y\in T(V)$ the \emph{braided commutator} as follows:
\begin{equation}\label{eqn:braidedcommutator}
[x,y]_c := \text{multiplication } \circ \left( \id - c \right)
\left( x \ot y \right).
\end{equation}

Fix a braided vector space of diagonal type $(V,c)$ and an ordered basis $X = \{x_1,\dots ,x_{\theta}\}$ of $V$ as in \eqref{ecuacion de trenzas}. Let $\xx$ be the corresponding vocabulary (the set of words with letters in $X$) and consider the lexicographical
order on $\xx$. We will identify the vector space $\kk \xx$ with $T(V)$. We shall consider two different gradings of the algebra $T(V)$. First,
the usual $\N_0$-grading $T(V) = \oplus_{n\geq 0}T^n(V)$. If we denote by $\ell$ the length of a word in $\xx$, then $T^n(V)=
\oplus_{x\in\xx, \, \ell(x) = n}\kk x$.
\label{paginatres}

Second, let $\alpha_1, \dots, \alpha_\theta$ be the canonical basis of $\zt$. Then $T(V)$ is $\zt$-graded, where the degree is given by $\deg x_i = \alpha_i$, $\otv$. Consider the bilinear form
$\chi: \zt\times \zt \to \kk^{\times}$ given by $\chi(\alpha_i, \alpha_j) = q_{ij}$,
$\otvz$. Then
\begin{equation}\label{braiding}
    c(u \ot v)= q_{u,v} v \ot u, \qquad u,v \in \xx,
\end{equation}
where $q_{u,v} = \chi(\deg u, \deg v)\in
\kk^{\times}$. The braided commutator satisfies a ``braided" derivation equation which gives place to a ``braided" Jacobi identity, namely
\begin{align}\label{idjac}
\left[\left[ u, v \right]_c, w \right]_c &= \left[u, \left[ v, w
\right]_c \right]_c
 - \chi( \alpha, \beta ) v \ \left[ u, w \right]_c + \chi( \beta, \gamma) \left[ u, w \right]_c \ v,
 \\
\label{der}
    \left[ u,v \ w \right]_c &= \left[ u,v \right]_c w + \chi( \alpha, \beta ) v \ \left[ u,w \right]_c,
\\ \label{der2} \left[ u \ v, w \right]_c &= \chi( \beta, \gamma ) \left[ u,w \right]_c \ v + u \ \left[ v,w \right]_c,
\end{align}
for any homogeneous $u,v,w \in T(V)$, of degrees $\alpha, \beta, \gamma \in \N^{\theta}$, respectively.

We denote by $\yd$ the category of Yetter-Drinfeld
modules over $H$, where $H$ is a Hopf algebra with bijective
antipode. Any $V\in \yd$ becomes a braided vector space \cite[Section 10.6]{M}.
If $H= \kk \Gamma$, where $\Gamma$ is a finite abelian group, then any
$V\in \yd$ is a braided vector space of diagonal type: if $V_{g} = \{v\in V \mid
\delta(v) = g\otimes v\}$, $V^{\chi} = \{v\in V \mid  g \cdot v =
\chi(g)v \text{ for all } g \in \Gamma\}$ and $V_{g}^{\chi} = V^{\chi} \cap V_{g}$, then $V
= \oplus_{g\in \Gamma, \chi\in \VGamma}V_{g}^{\chi}$.  In this setting the braiding is given
by
$$ c(x\otimes y) = \chi(g) y\otimes x, \qquad x\in V_{g}, \
g \in \Gamma, \ y\in V^{\chi}, \ \chi \in \VGamma.$$

Reciprocally, any braided vector space of diagonal type can be
realized as a Yetter-Drinfeld module over the group algebra of many abelian groups. For example let $(V,c)$ be a braided vector space of diagonal type. Call $\Gamma$ the
free abelian group of rank $\theta$, with basis $g_1, \ldots,
g_{\theta}$, and define the characters $\chi_1, \ldots,
\chi_{\theta}$ of $\Gamma$ by
    \[ \chi_j(g_i)=q_{ij}, \quad 1 \leq i,j \leq \theta. \]
We can consider $V$ as a Yetter-Drinfeld module over $\kk \Gamma$ for which $x_i \in V_{g_i}^{\chi_i}$.

\bigbreak Given $V\in \yd$, the tensor algebra $T(V)$ admits a
unique structure of graded braided Hopf algebra in $\yd$ such that the elements of 
$V$ are primitive. As in \cite{AS1}, we define the \emph{Nichols algebra} $\bB(V)$ associated to $V$ as the quotient of
$T(V)$ by the maximal element $I(V)$ of the family
$\bS$ of all the homogeneous two-sided ideals $I \subseteq \oplus_{n \geq 2} T(V)$
such that $I$ is a Yetter-Drinfeld submodule of $T(V)$ and a Hopf ideal: $\Delta(I) \subset I\ot T(V) + T(V)\ot I$.

The following proposition characterizes the Nichols algebra associated to $V$ in a very interesting way.

\begin{prop}\label{Prop:formabilineal} \cite[Prop. 1.2.3]{L3}, \cite[Prop. 2.10]{AS1}.
For each family of scalars $a_1, \ldots,a_{\theta} \in \kk^{\times}$, there exists a unique
bilinear form $( | ): T(V) \times T(V) \rightarrow \kk$ such that
$(1|1)=1$, and:
\begin{eqnarray}
    (x|yy') &=& (x_{(1)} | y) (x_{(2)} | y'), \quad \mbox{for all  } x,y,y' \in
    T(V);\label{bilinearprop2}
    \\ (xx'|y) &=& (x|y_{(1)}) (x'|y_{(2)}), \quad \mbox{for all  } x,x',y \in T(V), \label{bilinearprop3}
    \\ (x_i | x_j) &=& \delta_{ij}a_i, \quad \mbox{for all  } i,j. \label{bilinearprop1}
\end{eqnarray}
This form is symmetric and satisfies
\begin{equation}
    (x|y)=0, \quad \mbox{for all  } x \in T(V)_g, \ y \in T(V)_h, \ g,h \in \Gamma, \ g \neq h. \label{bilinearprop4}
\end{equation}

The radical of this form $\left\{ x \in T(V):
(x|y)=0, \ \forall y \in T(V) \right\}$ is $I(V)$, so $(\cdot|\cdot)$ induces a non degenerate bilinear form on $\bB(V)$
denoted by the same name. \qed
\end{prop}

In consequence, if $(V,c)$ is of diagonal type, then the ideal $I(V)$ is
$\zt$-homogeneous and $\bB(V)$ is $\zt$-graded, see \cite[Prop. 2.10]{AS1}.

\subsection{Lyndon words and PBW basis of braided graded Hopf algebras generated in degree zero and one}\label{subsection:pbw}

\

A word $u \in \xx$, $u\neq 1$, is \emph{Lyndon}
if $u$ is smaller than any of its proper ends; that is, for any decomposition $u=vw$, $v,w \in\xx - \left\{ 1 \right\}$, we have $u<w$. We denote by $L$ the set of Lyndon words, see \cite[Chapter 5]{Lo}

Note that $X \subset L$, and any Lyndon word begins by its smallest letter. They also satisfy the following properties.

\begin{enumerate}
    \item Let $u \in \xx-X$. Then $u$ is Lyndon if and only if for any
    decomposition $u=vw$, $v,w \in \xx - \left\{ 1 \right\}$, it satisfies $vw=u < wv$.
    \item If $v,w \in L, v<w$, then $vw \in L$.
    \item Let $u \in\xx-X$. Then $u \in L$ if and only if there exist
$v,w \in L$ with $v<w$ such that $u=vw$.
\end{enumerate}

\begin{defn}
Let $u \in L-X$. The \emph{ Shirshov decomposition } of $u$ is the decomposition $u=vw$, with $v,w \in L$
such that $w$ is the smallest end among those proper non-empty ends of $u$, see \cite{Lo}. Following \cite{He}, we denote it by $\Sh(u)=(v,w) \in L \times L$. It satisfies that $w$ is the longest end between the ends that are Lyndon words.

\end{defn}

Given $u,v,w \in L$ be such that $u=vw$, $u \neq 1$, then $\Sh(u)=(v, w)$ if and only if either $v \in X$, or else $\Sh(v)=(v_1, v_2)$ satisfies $w \leq
v_2$.

\bigbreak Lyndon Theorem says that any word $u \in \xx$ admits a unique decomposition $u=l_1l_2\dots  l_r$,
as a product of non-increasing Lyndon words: $l_i \in L$, $l_r \leq \dots \leq l_1$; see \cite[Thm. 5.1.5]{Lo}.  This is called the \emph{Lyndon
decomposition} of $u \in \xx$; we call \emph{Lyndon letters} of $u$ to any $l_i \in L$ appearing in such
decomposition.
\smallskip

We recall the endomorphism $\left[ - \right]_c$, see \cite{Kh}, defined inductively on $\kk \xx$ using Shirshov and Lyndon decomposition:
$$
\left[ u \right]_c := \begin{cases} u,& \text{if } u = 1
\text{ or }u \in X;\\
[\left[ v \right]_c, \left[ w \right]_c]_c,  & \text{if } u \in
L, \, \ell(u)>1 \text{ and }\Sh(u)=(v,w);\\
\left[ u_1 \right]_c \dots  \left[ u_t \right]_c,& \text{if } u
\in \xx-L \\ &\qquad
\text{ with Lyndon decomposition  }u=u_1\dots u_t.\\
\end{cases}
$$

\begin{defn}{\cite{Kh}.} The \emph{hyperletter} corresponding to
$l \in L$ is $\left[ l \right]_c$. A \emph{hyperword}
is a word in hyperletters, and a \emph{monotone hyperword} is a
hyperword $\left[u_1\right]_c^{k_1}\dots
\left[u_m\right]_c^{k_m}$ such that $u_1>\dots >u_m$.
\end{defn}

Note that for any $u \in L$, $\left[ u \right]_c$ is a homogeneous
polynomial with coefficients in the subring $\mathbb{Z} \left[q_{ij}\right]$
and $ \left[ u \right]_c\in u+ \mathbb{Z} \left[q_{ij}\right] \xx^{\ell(u)}_{>u}$.

\medbreak  The hyperletters inherit the order from the Lyndon
words; this induces in turn an ordering in the hyperwords (the
lexicographical order on the hyperletters). We describe now the braided commutator of hyperwords.

\begin{theorem}{\cite[Thm. 10]{R2}.} \label{Theo:corcheteentreLyndon} 
Let $m,n \in L$, with $m<n$. Then $\left[\left[m\right]_c, \left[n\right]_c \right]_c$ is a
$\mathbb{Z} \left[q_{ij}\right]$-linear combination of monotone
hyperwords $\left[l_1\right]_c \dots  \left[l_r\right]_c$, $l_i \in
L$, such that the hyperletters of those hyperwords satisfy $n>l_i \geq mn$.

Moreover, $\left[mn\right]_c$ appears in the expansion with
non-zero coefficient, and for any hyperword of this
decomposition, $\deg (l_1\dots l_r)= \deg(mn)$. \qed
\end{theorem}

The coproduct of $T(V)$ can be described also in the basis of hyperwords.

\begin{lema}\label{Lemma:coproductPBWelements} \cite{R2}.
Let $u \in \xx$, and $u= u_1\dots u_r v^m, \ v, u_i \in L, v<u_r
\leq \dots  \leq u_1$ be the Lyndon decomposition of $u$. Then
\begin{eqnarray*}
        \Delta \left(\left[ u \right]_c\right) &=& 1 \ot \left[ u \right]_c+ \sum ^{m}_{i=0} \binom{ m }{ i } _{q_{v,v}} \left[u_1\right]_c\dots  \left[u_r\right]_c \left[ v \right]_c ^i \ot \left[ v \right]_c^{m-i}
        \\ && + \sum_{ \substack{ l_1\geq \dots  \geq l_p >v, \  l_i \in L \\ 0\leq j \leq m } } x_{l_1,\dots ,l_p}^{(j)}
        \ot \left[l_1\right]_c\dots
        \left[l_p\right]_c\left[v\right]_c^j;
\end{eqnarray*}
where each $x_{l_1,\dots ,l_p}^{(j)}$ is $\zt$-homogeneous, $\deg(x_{l_1,\dots ,l_p}^{(j)} \, l_1\dots  l_p v^j)= \deg(u)$.
\qed
\end{lema}

We have then the following result from \cite{R2}.

\begin{lema}\label{Lemma:subalgebrasWl}
For each $l\in L$ call $W_l$ the subspace of $T(V)$ generated by
\begin{equation}\label{elementsWl}
    [l_1]_c [l_2]_c\cdots [l_k]_c, \quad k \in \N_0, \, l_i \in L, \, l_1 \geq \ldots \geq l_k \geq l.
\end{equation}
Then $W_l$ is a right coideal subalgebra of $T(V)$.
\end{lema}
\bdem
It follows from Theorem \ref{Theo:corcheteentreLyndon} and Lemma \ref{Lemma:coproductPBWelements}.
\edem

We consider as in \cite{U} and \cite{Kh} another order in $\xx$. Given $u,v \in \xx$, we say that $u \succ v$ if
and only if either $\ell(u)<\ell(v)$, or else $\ell(u)=\ell(v)$ and $u>v$ for the lexicographical order.
We call $\succ$ the \emph{deg-lex order}, which is a total order. The empty word 1 is the maximal element for
$\succ$, and this order is invariant by right and left multiplication.

\medskip

Let $I$ be a proper ideal of $T(V)$, and set $R=T(V)/I$.
Let $\pi: T(V) \rightarrow R$ be the canonical projection. Let us
consider the subset of $\xx$:
    \[G_I:= \left\{ u \in \xx: u \notin \\ \kk \xx_{\succ u}+I  \right\}. \]
Such set satisfies:
\begin{enumerate}
    \item[(a)] If $u \in G_I$ and $u=vw$, then $v,w \in G_I$.
    \item[(b)] Any $u \in G_I$ factorizes uniquely as a non-increasing product of Lyndon words in $G_I$.
\end{enumerate}

\begin{prop}{\cite{Kh,R2}.}\label{firstPBWbasis}
The set $\pi(G_I)$ is a basis of $R$. \qed
\end{prop}

\noindent In what follows, we assume that $I$ is a Hopf ideal. Consider now
\begin{equation}\label{setsi}
S_I:= G_I \cap L.
\end{equation}
We then define the \emph{height} function $h_I: S_I \rightarrow \left\{2,3,\dots
\right\}\cup \left\{ \infty \right\}$ by
\begin{equation}\label{defheight}
    h_I(u):= \min \left\{ t \in \N : u^t  \in \kk \xx_{\succ u^t} + I \right\}.
\end{equation}
One can find a PBW-basis by hyperwords of the quotient $R$ of $T(V)$ using the set $S_I$ and the height previously defined.

\begin{theorem}{\cite{Kh}.}\label{basePBW}  The following set is a PBW-basis of $R=T(V)/I$:

$ \{ [u_1]_c^{n_1} \cdots [u_k]_c^{n_k}: \, k \in \N_0, \, u_1 > \ldots > u_k \in S_I, \, 0 \leq n_i <h_I(u_i) \} $. \qed
\end{theorem}
Proofs are in \cite{Kh}, where the next consequences are also considered.

\begin{prop}\label{altf}
For any $v \in S_I$ such that $h_I(v)< \infty$, $q_{v,v}$ is a root of unity, whose order coincides with $h_I(v)$. \qed
\end{prop}

\begin{cor}\label{cor:primero}
A word $u$ does not belong to $G_I$ if and only if the associated
hyperletter $\left[u\right]_c$  is a linear combination,
modulo $I$, of hyperwords $\left[ w \right]_c$, $w \succ u$, where
all the hyperwords have their hyperletters in $S_I$.

Moreover, if $h_I(v):= h < \infty$, then $\left[ v \right]^{h}$ is a linear
combination of hyperwords $\left[ w \right]_c$, $w \succ v^h$.
\qed
\end{cor}

\bigbreak

\section{Root systems and coideal subalgebras}\label{section:rootsystems-coidealsubalgebras}

In this section we recall the definition of Weyl groupoid and the associated generalized root system given in \cite{CH1} and \cite{HY}. We recall also some properties of these objects that we shall use in the subsequent sections, and the relation with Nichols algebras of diagonal type. After that, we describe convex orders for subsets of the root systems as a generalization of Papi's results in \cite{P} for Weyl groups. We consider then a family of coideal subalgebras of a Nichols algebra of diagonal type with finite root system in order to prove that the ordering on the Lyndon words of a PBW basis as in Section \ref{subsection:pbw} is convex. For the proof of the convexity we use the characterization of coideal subalgebras given in \cite{HS}.

\subsection{Weyl groupoid and root systems}\label{subsection:weylgroupoid}

The notation used here is the same as in \cite{CH1}.

Fix a non-empty set $\cX$, a non-empty finite set $I$ and call $\{\alpha_i \}_{i \in I}$ the canonical basis of $\Z^I$. For each $i \in I$ consider a map $r_i: \cX \rightarrow \cX$, and for each $X \in \cX$ a generalized Cartan matrix $A^X= (a^X_{ij})_{i,j \in I}$ in the sense of \cite{K}.

\begin{defn}{\cite{HY,CH1}}
The quadruple $\cC:= \cC(I, \cX, (r_i)_{i \in I}, (A^X)_{X \in \cC})$ is a \emph{Cartan scheme} if
\begin{itemize}
  \item for all $i \in I$, $r_i^2=id$, and
  \item for all $X \in \cX$ and $i,j \in I$: $a^X_{ij}=a^{r_i(X)}_{ij}$.
\end{itemize}
For each $i \in I$ and $X \in \cX$ denote by $s_i^X$ the automorphism of $\Z^I$ given by $$ s_i^X(\alpha_j)=\alpha_j-a_{ij}^X\alpha_i, \qquad j \in I. $$
The \emph{Weyl groupoid} of $\cC$ is the groupoid $\cW(\cC)$ whose set of objects is $\cX$ and whose morphisms are generated by $s_i^X$, where we consider $s_i^X \in \Hom(X, r_i(X))$, $i \in I$, $X \in \cX$.
\end{defn}

In general we shall denote $\cW(\cC)$ simply by $\cW$, and for any $X \in \cX$:
\begin{align}\label{defHom}
    \Hom(\cW,X) & := \cup_{Y \in \cX} \Hom(Y,X),
    \\ \label{defrealroot} \de^{X \ re} &:= \{ w(\alpha_i): \ i \in I, \ w \in \Hom(\cW,X) \}.
\end{align}
$\de^{X \ re}$ is the set of \emph{real roots} of $X$. Each $w \in \Hom(\cW,X_1)$ can be written as a product $s_{i_1}^{X_1}s_{i_2}^{X_2} \cdots s_{i_m}^{X_m}$, where $X_j=r_{i_{j-1}} \cdots r_{i_1}(X_1)$, $i \geq 2$. We denote $w= \id_{X_1} s_{i_1} \cdots s_{i_m}$: this means that $w \in \Hom(\cW,X_1)$, because the elements $X_j \in \cX$ are univocally determined. The \emph{length} of $w$ is defined by
$$ \ell(w)= \min \{ n \in \N_0: \ \exists i_1, \ldots, i_n \in I \mbox{ such that }w=\id_X s_{i_1} \cdots s_{i_n} \}. $$
In what follows we will assume that $\cW$ is a connected groupoid: $$\Hom(Y,X) \neq \emptyset, \quad \forall X,Y \in \cX.$$

\begin{defn}{\cite{HY,CH1}}
Fix a Cartan scheme $\cC$, and for each $X \in \cX$ a set $\de^X \subset \Z^I$. $\cR:= \cR(\cC, (\de^X)_{X \in \cX} )$ is a \emph{root system of type} $\cC$ if
\begin{enumerate}
  \item for all $X \in \cX$, $\de^X= (\de^X \cap \N_0^I) \cup  -(\de^X \cap \N_0^I)$,
  \item for each $i \in I$ and each $X \in \cX$, $\de^X \cap \Z \alpha_i= \{\pm \alpha_i \}$,
  \item for each $i \in I$ and each $X \in \cX$, $s_i^X(\de^X)=\de^{r_i(X)}$,
  \item if $m_{ij}^X:= |\de^X \cap (\N_0\alpha_i+\N_0 \alpha_j)|$, then $(r_ir_j)^{m_{ij}^X}(X)=(X)$ for all $i \neq j \in I$ and all $X\in \cX$.
\end{enumerate}
We call $\de^X_+:= \de^x \subset \N_0^I$ the set of \emph{positive roots}, and $ \de^X_-:= - \de^X_+$ the set of \emph{negative roots}.
\end{defn}
By (3) we have that $w(\de^X)= \de^Y$ for any $w \in \Hom(Y,X)$.

We say that $\cR$ is \emph{finite} if $\de^X$ is finite for some $X\in \cX$. By \cite[Lemma 2.11]{CH1}, this is equivalent to the fact that the sets $\de^X$ are finite, for all $X \in \cX$, and that $\cW$ is finite.

The following result plays a fundamental role for our purposes in the next subsection.

\begin{theorem}{\cite[Thm. 2.10]{CH2}}  \label{theorem:decomposition-of-roots}
Let $\alpha \in \Delta^X_+ \setminus \{ \alpha_i: i=1, \ldots \theta\}$. There exist $\beta, \gamma \in \Delta^X_+$ such that $\alpha=\beta+ \gamma$. \qed
\end{theorem}

Now we recall some results involving real roots and the length of the elements in $\cW$.

\begin{lema}{\cite[Cor. 3]{HY}}  \label{Lemma:lengthHY}
Let $m \in \N$, $X, Y \in \cX$ and $i_1, \ldots, i_m,j \in I$. Call $w=\id_X s_{i_1} \cdots s_{i_m} \in \Hom(Y,X)$, and assume that $\ell (w)=m$. Then,
\begin{itemize}
  \item $\ell (w s_j)=m+1$ if and only if $w(\alpha_j) \in \Delta^X_+$,
  \item $\ell (w s_j)=m-1$ if and only if $w(\alpha_j) \in \Delta^X_-$.
\end{itemize}
\qed
\end{lema}

\begin{prop}{\cite[Prop. 2.12]{CH1}}  \label{Prop:maxlengthCH}
For each $w=\id_X s_{i_1} \cdots s_{i_m}$ such that $\ell(w)=m$, the roots $\beta_j=s_{i_1} \cdots s_{i_{j-1}}(\alpha_{i_j}) \in \de^X$ are positive and pairwise different. If $w$ is an element of maximal length and $\cR$ is finite, then $\{ \beta_j \}= \Delta^X_+$. In consequence, all the roots are real: i.e., for each $\alpha \in \Delta^X_+$ there exist $i_1, \ldots, i_k, j \in I$ such that
$ \alpha = s_{i_k} \cdots s_{i_1}(x_j)$. \qed
\end{prop}

As in \cite{HS}, consider for $X \in \cX$, $m \in \N$ and $(i_1, \ldots, i_m) \in I^m$ the sets:
\begin{align}
\label{defnLambda1} \Lambda^X(i_1, \ldots, i_m) &:= \{ \beta_k:= \id_X s_{i_1} \cdots s_{i_k-1}(\alpha_{i_k}): \ 1 \leq k \leq m \} \subset \de^X,
\\ \label{defnLambda2} \Lambda_+^X(i_1, \ldots, i_m) &:= \{\beta \in \de^X_+: \ |\{k \in \{1, \ldots, m \}: \beta= \pm \beta_k \}| \mbox{ is odd} \}.
\end{align}
By \cite[Prop. 1.9]{HS}, given other elements $j_1, \ldots, j_n \in I$, we have
$$\Lambda_+^X(i_1, \ldots, i_m) = \Lambda_+^X(j_1, \ldots, j_n) \ \Leftrightarrow \ \id_x s_{i_1} \cdots s_{i_m} = \id_x s_{j_1} \cdots s_{j_n},$$and moreover,
\begin{equation}\label{cardinal-Lambda-set}
    |\Lambda_+^X(i_1, \ldots, i_m)| = \ell (\id_x s_{i_1} \cdots s_{i_m}).
\end{equation}
In this way, if $w=\id_X s_{i_1} \cdots s_{i_m}$ is any expression of $w \in \Hom(\cW,X)$, we can define $\Lambda_+^X(w):=\Lambda_+^X(i_1, \ldots, i_m)$.
\medskip

\subsection{Convex orders on root systems}\label{subsetion:convexorder}

Now we characterize convex orders on subsets of root systems of finite Weyl groupoids, extending the results given in \cite{P} for Weyl groups.

\begin{defn}
Consider a root system $\Delta^X_+$ with a fixed total order $<$. We say that it is

$\bullet$ \emph{convex} if for each $\alpha, \beta \in \Delta^+$ such that $\alpha <  \beta$ and $\alpha+\beta \in \Delta^+$, then
 $$ \alpha < \alpha+\beta <\beta. $$

$\bullet$ \emph{sub-convex} if for each $\alpha, \beta \in \Delta^+$ such that $\alpha <  \beta$ and $\alpha+\beta \in \Delta^+$, then
 $$ \alpha < \alpha+\beta. $$

$\bullet$ \emph{strongly convex} if for each ordered subset $\alpha_1 \leq  \ldots \leq \alpha_k$ of $\Delta^{+}$ such that $\alpha := \sum \alpha_i \in  \Delta^+$ then $\alpha_1 < \alpha < \alpha_k$.
\end{defn}

\begin{defn} Let $L= \{ \beta_1, \ldots, \beta_m \}$ be an ordered subset of $\Delta^X_+$. We say that $L$ is \textit{associated to} $w \in \Hom(\cW, X)$ if there exists a reduced expression $w=\id_X s_{i_1} \cdots s_{i_m}$ such that
$$\beta_j=s_{i_1} \cdots s_{i_{j-1}}(\alpha_{i_j}), \qquad \forall 1 \leq j \leq m.$$
Compare this with \cite{P}. For any $w \in \Hom(Y, X)$ define
$$ R_{w}:= \{ \alpha \in \Delta^X_+: w^{-1}(\alpha) \in \Delta^Y_- \} .$$
\end{defn}

Now we generalize some results about Weyl groups to the context of Weyl groupoids. First we consider the analogue of a result in \cite{B}.

\begin{prop}\label{Prop:Bourbaki}
For any ordered set $L$ associated to $w$, we have $L=R_{w}$. In consequence, $|R_{w}|= \ell(w)$ and two ordered sets associated to the same $w$ differ at most by the ordering.
\end{prop}
\bp
Note that for any $\beta_j=s_{i_1} \cdots s_{i_{j-1}}(\alpha_{i_j})$, $w^{-1}(\beta_j)= -s_{i_m} \cdots s_{i_{j+1}} (\alpha_j)$. $s_{i_m} \cdots s_{i_{j+1}} s_{i_j}$ is a reduced expression because it is contained in a reduced expression, so we have $w^{-1}(\beta_j) \in \Delta^Y_-$ by Lemma \ref{Lemma:lengthHY}. Therefore $L \subseteq R_{w}$.

Reciprocally, let $\alpha \in R_{w}$. As $w^{-1}(\alpha) \in \de^Y_-$ and $s_{i_1} \cdots s_{i_m} (w^{-1}(\alpha))= \alpha \in \de^X_+$, consider the greatest $j$ such that $s_{i_j} \cdots s_{i_m} w^{-1}(\alpha)$ is positive. Therefore $s_{i_{j+1}} \cdots s_{i_m} w^{-1}(\alpha)$ is negative, so $s_{i_j} \cdots s_{i_m} w^{-1}(\alpha) = \alpha_{i_j}$, and then $\alpha_{i_j}= s_{i_j} \cdots s_{i_m} w^{-1}(\alpha)$; that is, $\alpha=s_{i_1} \cdots s_{i_{j-1}}(\alpha_j) \in L$.
\ep

Second, we relate our sets $R_w$ with the ones in \cite{HS}, see \eqref{defnLambda2}. Although the sets are equal, our definition is more comfortable to prove statements about convexity.

\begin{lema}\label{Lemma:R-Lambda-iguales}
For each $w \in \Hom(\cW,X)$, $R_w=\Lambda_+^X(w).$
\end{lema}
\bp
Fix a reduced expression $w=\id_X s_{i_1} \cdots s_{i_m}$, so $\beta_j=s_{i_1} \cdots s_{i_{j-1}}(\alpha_{i_j})$ is a positive root, and $\alpha \in \de^X_+$ is equal to $\pm \beta_j$ if and only if $\alpha=\beta_j$. Therefore $\Lambda_+^X(w)=L$.
\ep

Now we extend another result from \cite{P}. Note that condition (a) in our result is weaker than the one in \cite{P}, but the proof is very similar. This weaker condition shall simplify some proofs in what follows.

\begin{theorem}\label{Thm:Papigeneralizado}
Let $L$ be an ordered subset of $\Delta^X_+$. There exists $w \in \Hom(\cW,X)$ such that $L$ is associated to $w$ if and only if the following conditions are satisfied:
\begin{enumerate}
  \item[(a)] For each pair $\lambda < \mu \in L$ such that $\lambda+\mu \in \Delta^X_+$, then $\lambda + \mu \in L$ and $\lambda < \lambda + \mu$.
  \item[(b)] If $\lambda + \mu \in L$ and $\lambda, \mu \in \Delta^X_+$, then at least one of them belongs to $L$ and precedes $\lambda + \mu$.
\end{enumerate}
\end{theorem}
\bp
Assume that $L$ is associated to $w=\id_X s_{i_1} \cdots s_{i_m}$ for some $w \in \Hom(Y,X)$. If $\lambda=s_{i_1} \cdots s_{i_{k-1}} (\alpha_{i_k})$ and $\mu=s_{i_1} \cdots s_{i_{j-1}} (\alpha_{i_j})$ are such that $1 \leq k <j \leq m$ and $\lambda+\mu \in  \Delta^X_+$, we have $\lambda + \mu \in L=R_w$, because $$ w^{-1}(\lambda+\mu)=w^{-1}(\lambda) + w^{-1}(\mu) \in \Delta^Y_- .$$
Suppose that $\lambda+ \mu < \lambda$. Then $\lambda+\mu= s_{i_1} \cdots s_{i_{l-1}} (\alpha_{i_l})$ for some $1 \leq l <k$, so $s_{i_l} \cdots s_{i_1} (\lambda+\mu)=-\alpha_l \in \Delta^{r_{i_l} \cdots r_{i_1}(X)}_-$. But as $l<k<j$, we have $$s_{i_l} \cdots s_{i_1}(\lambda), \, s_{i_l} \cdots s_{i_1}(\mu) \in \Delta^{r_{i_l} \cdots r_{i_1}(X)}_- ,$$
which is a contradiction. Therefore $\lambda < \lambda + \mu$, and $L$ satisfies (a).

For (b), suppose that $\lambda+\mu \in L$, but $\lambda, \mu \notin L$: $w^{-1}(\lambda), w^{-1}(\mu) \in \Delta^Y_+$, so $w^{-1}(\lambda+\mu)$ is positive, which is a contradiction to the fact that $\lambda+\mu \in R_w$. If both $\lambda, \mu \in L$, a similar proof to (a) gives that one of them precedes $\lambda+\mu$. In consequence, suppose that $\lambda \in L, \mu \notin L$ and $\lambda+ \mu < \lambda$. If $l<k$ is such that $\lambda+\mu= s_{i_1} \cdots s_{i_{l-1}} (\alpha_{i_l})$, we have $s_{i_l} \cdots s_{i_1}(\lambda) \in \Delta^+$ and
$$ s_{i_l} \cdots s_{i_1}(\lambda)+ s_{i_l} \cdots s_{i_1}(\mu) = s_{i_l} \cdots s_{i_1}(\lambda + \mu)=-\alpha_l \in \Delta^{r_{i_l} \cdots r_{i_1}(X)}_-,$$
so $s_{i_l} \cdots s_{i_1}(\mu) \in \Delta^{r_{i_l} \cdots r_{i_1}(X)}_-$, and then $\mu \in R_{\id_X s_{i_1} \cdots s_{i_l}} \subset R_{\id_X s_{i_1} \cdots s_{i_m}}=L$, a contradiction.

\medskip

Reciprocally, we prove that an ordered set $L$ satisfying (a) and (b) is associated to some $w$ by induction on $m:=|L|$. If $m=1$, let $\alpha \in L$. If we suppose that $\alpha$ is not simple, by Theorem \ref{theorem:decomposition-of-roots}, $\alpha=\beta+\gamma$ for some positive roots $\beta, \gamma$, and by condition (b) one of them belongs to $L$, so $m \geq 2$, which is a contradiction. Therefore $L= \{ \alpha_j \} = R_{s_j}$ for some $1 \leq j \leq \theta$.

Now assume $m > 1$ and call $\beta_1 < \ldots < \beta_m$ the elements of $L$. Notice that $L'=\{ \beta_1, \ldots, \beta_{m-1} \}$ verifies conditions (a) and (b), so by inductive hypothesis there exists a reduced expression $v=\id_X s_{i_1} \cdots s_{i_{m-1}}$ such that
$$ \beta_1= \alpha_{i_1}, \qquad \beta_j= s_{i_1} \cdots s_{i_{j-1}}(\alpha_{i_j}), \quad j=2,\ldots,m-1. $$
Let $Z=r_{i_{m-1}} \cdots r_{i_1}(X)$. Then $v^{-1}(\beta_m) \in \de^Z_+$ because $\beta_m \notin L'=R_v$. Suppose that $v^{-1}(\beta_m)$ is not simple. Then there exist $\alpha, \beta \in \de^Z_+$ such that $\alpha+\beta = v^{-1}(\beta_m)$; i.e. $\beta_m= \alpha'+\beta'$, where $\alpha'=v(\alpha),\ \beta'=v(\beta) \in \de^X$. Therefore $\alpha' \in \de^X_+$ or $\beta' \in \de^X_+$. On the other hand, if both are positive then one of them is $\beta_k$ for some $k < m$; assume $\alpha'=\beta_k$, but then $\alpha=v^{-1} (\beta_k) \in \de^Z_-$, a contradiction. In consequence, we can consider $\alpha' \in \de^X_+$ and $\beta' \in \de^X_-$. For this case, $\alpha' \notin R_v=L'$ and $-\beta' \in R_v=L' \subset L$. As $\alpha'= \beta_m+(-\beta')$, hypothesis (a) on the set $L$ implies that $\alpha' \in L$, so $\alpha'=\beta_m \in L -L'$, a contradiction. Therefore, $v^{-1}(\beta_m) = \alpha_{i_m}$ for some $i_m \in I$, $w= v s_{i_m} \in \Hom (r_{i_m}(Z), X)$ is a reduced expression by Lemma \ref{Lemma:lengthHY}, and $L=R_w$.
\ep

\begin{theorem}\label{Theorem:equivalenciasconvexo}
Given an order on $\Delta^X_+$, the following statements are equivalent:
\begin{enumerate}
  \item the order is associated with a reduced expression of the longest element,
  \item the order is strongly convex,
  \item the order is convex.
\end{enumerate}
\end{theorem}
\bp
(1) $\Rightarrow$ (2). Let $\omega= \id_X s_{i_1} \cdots s_{i_m}$ be an element of maximal length in $\Hom(\cW,X)$. By Proposition \ref{Prop:maxlengthCH}, $m= |\Delta^X_+|$ and
$$ \beta_k := s_{i_1} \cdots s_{i_{k-1}} (\alpha_{i_k}), \quad k=1, \ldots, m, $$
are all different, so $\{ \beta_k\}= \Delta^X_+$. In consequence, it induces an order on $\Delta^X_+$: $$\beta_1 < \cdots < \beta_m.$$

To prove that this order is strongly convex, consider $\beta, \beta_{k_1}, \ldots, \beta_{k_l} \in \Delta^X_+$ such that $k_1 < \cdots <k_l$ and $\beta= \beta_{k_1}+ \cdots +\beta_{k_l}$. Suppose that $\beta=\beta_k$ with $k<k_1$. Then $u=\id_X s_{i_1} \cdots s_{i_k}$ satisfies $\ell(u)=k$, $\beta \in R_u$ but $\beta_{k_j} \notin R_u$ for all $j=1,...,l$, which is a contradiction because $u(\beta) \in \de^{r_{i_k} \cdots r_{i_1}X}_-$ should be the sum of the positive roots $u(\beta_j)$. We obtain a similar contradiction if we assume $k>k_l$. Therefore $k_1 <k<k_l$.

\smallskip

(2) $\Rightarrow$ (3) is clear.
\smallskip

(3) $\Rightarrow$ (1). Assume that a given order on $\Delta^X_+$ is convex; i.e. it satisfies trivially condition (a) of Theorem \ref{Thm:Papigeneralizado} because we consider $L=\Delta^X_+$. Therefore it satisfies also condition (b) by the convexity, so the order is associated to some $w$. As $\ell(w)=|\Delta^X_+|$ by Proposition \ref{Prop:Bourbaki}, it should be the element of maximal length.
\ep

\subsection{Coideal subalgebras and convex orders for PBW bases}\label{subsection:coidealsubalgebras}

Now we recall a description of coideal subalgebras of Nichols algebras with finite root system given in \cite{HS}. We will use this result to prove that the lexicographical order on the PBW generators of Kharchenko's basis is convex. Before to prove it, we recall the results about the Weyl groupoid attached to a braided vector space of diagonal type. Given a braided vector space $(V,c)$ of diagonal type, fix a basis $\{ x_1, \ldots,x_{\theta}\}$ and   scalars $q_{ij} \in \kk^\times$ as in \eqref{ecuacion de trenzas},  and the bilinear form as in \eqref{braiding}. We set as in \cite{H1},
$\Delta^V_+ $ the set of degrees of a PBW basis of $\bB(V)$, counted with their multiplicities. Such set does not depend on the PBW basis, as it is remarked in \cite{H1} and proved in \cite{AA}.

In what follows, \textbf{we fix a braided vector space $(V,c)$ of diagonal type and assume that the root system $\de^V_+$ is finite.}  In such case we can attach a Cartan scheme $\cC$, a Weyl groupoid $\cW$ and the corresponding root system $\cR$, see \cite[Thms. 6.2, 6.9]{HS} and the references therein, which coincides with the Weyl groupoid defined in \cite{H1} for braided vector space of diagonal type. Such Weyl groupoid can be built as follows, see \cite{AA}.  Set $\cX$ the set of ordered bases of $\zt$, and for each $F= \left\{ f_1, \ldots,f_{\theta} \right\} \in \cX$, set $\widetilde{q}_{ij}= \chi(f_i,f_j)$. Define for each $1 \leq i \neq j \leq \theta$,
\begin{equation}
m_{ij}(F):= \min \left\{ n \in \mathbb{N}_0: (n+1)_{\widetilde{q}_{ii}} (1-\widetilde{q}_{ii}^n \widetilde{q}_{ij} \widetilde{q}_{ji} )=0 \right\}, \label{mij}
\end{equation}
and set $s_{i,F} \in \Aut(\mathbb{Z}^{\theta})$ such that $s_{i,F}(f_j)=f_j+m_{ij}(F)f_i$. Here $m_{ii}=-2$.

Note that $\cG=\Aut(\mathbb{Z}^{\theta}) \times \cX$ is a groupoid whose set of objets is $\cX$ and whose morphisms are:
\begin{center}
$x \stackrel{(g,x)}{\longrightarrow} g(x)$.
\end{center}
The \emph{Weyl Groupoid} $W(\chi)$ of $\chi$ is the least subgroupoid of $\cG$ such that
\begin{itemize}
    \item $(id,E) \in W(\chi)$,
    \item if $(id,F) \in W(\chi)$ and $s_{i,F}$ is defined, then $(s_{i,F},F) \in W(\chi)$.
\end{itemize}
\bigskip

Recall that the (right) \emph{Duflo order} on $\Hom(\cW, X)$ is defined as follows: if $x \in  \Hom(Y,X)$ and $y \in \Hom(Z,Y)$, then $x \leq_D xy$ iff $\ell(xy)= \ell(x) + \ell(y)$; see \cite[Defn. 1.11]{HS}. By \cite[Thm. 1.13]{HS}, given $v,w \in \Hom(\cW,X)$ we have $v \leq_D w$ if and only if $\Lambda^X_+(v) \subset \Lambda^X_+(w)$.

\begin{obs}\label{Rem:maximallength}
Let $w_1 \leq_D w_2 \leq_D \cdots \leq_D w_k$ be a maximal chain in $\Hom(\cW,X)$. Then there exist a reduced expression $\id_X s_{i_1} \cdots s_{i_k}$ for some $i_1, \ldots, i_k \in I$ such that $w_j=\id_X s_{i_1} \cdots s_{i_j}$, for each $1 \leq j \leq k$.

In particular, a chain $w_1 \leq_D w_2 \leq_D \cdots \leq_D w_k$ has maximal length iff it is associated to a reduced expression of the longest element in $\Hom(\cW,X)$, and in consequence $k=|\de^X_+|$.
\end{obs}

We recall now some results from \cite{HS} about the classification of coideal subalgebras of $\bB(V)$. As in loc. cit., we
denote by $\cK(V)$ the set of all the $\N_0^\theta$-graded left coideal subalgebras of $\bB(V)$. We rewrite these results in
the context of diagonal braidings (in \cite{HS} the authors work in a more general context).

First results about the classification of coideal subalgebras were obtained in \cite{Kh2,KhL,Po} for quantized enveloping algebras
$U_q(\mathfrak{g})$ of type $A_n$, $B_n$ and $G_2$, respectively, where it was proved that coideal subalgebras admit a PBW basis and these subalgebras were classified.

\smallskip

Given $n=(n_1, \ldots , n_{\theta}) \in \N_0^\theta$, we set $X^{n}=X_1^{n_1} \cdots
X_{\theta}^{n_{\theta}}$ in $\kk [[x_1, \ldots, x_{\theta}]]$. We also set
\begin{eqnarray*}
\bq _h(\mathrm{t}) := \frac{\mathrm{t}^h-1}{\mathrm{t}-1} \in \kk
[\mathrm{t}], \quad h \in \N; \qquad \bq _{\infty}(\mathrm{t}):=
\frac{1}{1-\mathrm{t}}= \sum_{s=0}^{\infty} \mathrm{t}^s \in \kk
[[\mathrm{t}]].
\end{eqnarray*}
For each $\N^\theta_0$-graded $\kk$-vector space $W=\oplus_{\alpha \in \N^\theta_0} W_{\alpha}$, we denote its Hilbert series by
$$ \cH_W := \sum_{\alpha \in \N^\theta_0} (\dim W_{\alpha}) X^\alpha \in  \kk [[x_1, \ldots, x_{\theta}]]  .  $$

For any $\alpha \in \N_0^\theta$, we set $q_\alpha = \chi(\alpha, \alpha)$, where $\chi$ is the bicharacter over $\Z^\theta$
as in \eqref{braiding}, and $N_\alpha= \ord q_\alpha$, where $N_\alpha= \infty$ if $q_\alpha$ is not a root of unity.

\begin{theorem}\cite{HS}\label{Thm:HSadaptado}
For each $w \in \Hom(\cW,V)$ there exists a unique left coideal subalgebra $F(w) \in \cK(V)$ such that its Hilbert series is
\begin{equation}\label{HilbertseriesF}
    \cH_{F(w)} = \prod_{\beta \in \Lambda^V_+(w)} \bq_{N_\beta}(X^{\beta}).
\end{equation}

Moreover, the correspondence $w \mapsto F(w)$ gives an order preserving and order reflecting bijection between $\Hom(\cW,V)$ and $\cK(V)$, where we consider the Duflo order over $\Hom(\cW,V)$ and the inclusion order over $\cK(V)$; i.e.
$$ w_1 \leq_D w_2 \ \Leftrightarrow \ F(w_1) \subset F(w_2) . $$
\end{theorem}
\bp
Note that in \cite{HS} the authors classify right coideal subalgebras, but that $E$ is a right coideal subalgebras if and
only if $\cS(E)$ is a left coideal subalgebra, where $\cS$ denotes the antipode of $\bB(V)$. Moreover, if they are
$\N^{\theta}$-graded, then $\cH_E=\cH_{\cS(E)}$, because $\cS$ is $\N^{\theta}$-graded, and the order given by inclusion
on the family of left coideal subalgebras corresponds with the one on the family of right coideal subalgebras because $\cS$
is bijective. In this context we define $F(w)=\cS(E^V(w))$, where $E^V(w)$ is as in \cite[Thm. 6.12]{HS}.

By \cite[Lemma 6.11]{HS}, we have an isomorphism of $\N_0^\theta$-graded spaces $F(w) \cong \ot_{\beta \in \Lambda^V_+(w)}
\bB(V_\beta)$, where $V_\beta$ corresponds to $N_{\beta}$ of \cite[Defn. 6.5]{HS}. In this way $V_\beta$ is a 1-dimensional
braided vector space of diagonal type generated by a non-zero vector $v_\beta$, such that $c(v_\beta \ot v_\beta)=q_{\beta}
\ v_\beta\ot v_\beta$. Therefore, $\cH_{\bB(V_\beta)} = \bq_{N_\beta}(X^{\beta})$, and equation \eqref{HilbertseriesF}
follows.

The uniqueness of a coideal subalgebra with a given Hilbert series follows from \cite[Lemma 6.4]{HS}. The map $\Hom(\cW,V)
\rightarrow \cK(V)$ is bijective and preserves the order in both directions by \cite[Thms. 6.12, 6.15]{HS} (note that we can
apply these Theorems because we assume that $V$ has diagonal braiding and  $\de^V_+$ is finite).
\ep

\bigskip

Consider the PBW basis of Lyndon words given in Theorem \ref{basePBW} for the fixed basis $\{x_1, \ldots, x_\theta \}$ of
$V$. We assume that $\de^V_+$ is finite, so all the roots are real and have multiplicity one. In this way, we can label
the PBW generators by the elements $\beta \in \Delta^V_+$: they are $x_{\beta}= [l_{\beta}]_c$ for some Lyndon word
$l_\beta$ of degree $\beta$. It induces a total order on the roots: if $l_{\beta_1} < l_{\beta_2} < \cdots < l_{\beta_M}$
are ordered lexicographically, we consider $\beta_1 < \beta_2 < \cdots < \beta_M$, where $M= |\de^V_+|$ and in particular
$l_{\beta_1}=x_1$, $l_{\beta_M}=x_\theta$. Call $B$ the basis of $\bB(V)$ consisting of hyperwords as above.

Let $\pi: T(V) \rightarrow \bB(V)=T(V)/I(V)$ be the canonical projection. Recall the definition of the coideal subalgebras
$W_{l_\beta}$ in Lemma \ref{Lemma:subalgebrasWl}, and call
$$ W_\beta:= \pi(W_{l_{\beta}}), \qquad \beta \in \de^V_+. $$

\begin{obs}\label{Rem:Walpha}
$W_\beta$ is a left coideal subalgebra of $\bB(V)$, because $\pi$ is a morphism of braided Hopf algebras and $W_{l_\beta}$
is a left coideal subalgebra of $T(V)$. Also $W_{\beta_j} \subseteq W_{\beta_i}$ if $i < j$ and
$$ W_{\beta_1}=\bB(V), \quad W_{\beta_M}= \kk \langle x_\theta \rangle. $$
\end{obs}

\begin{lema}\label{Lemma:Walphacadenamaximal} With the notation above, $x_{\beta_i} \notin W_{\beta_j}$ if $i < j$.

In consequence, $\bB(V)= W_{\beta_1} \supsetneq W_{\beta_2} \supsetneq \cdots \supsetneq W_{\beta_M}$.
\end{lema}

\bp
Suppose that $x_{\beta_i} \in W_{\beta_j}$ with $i < j$. Then $x_{\beta_i} \in G_{I(V)}$ is a linear combination of
hyperwords greater or equal that $x_{\beta_j}$ in $\bB(V)$, which is a contradiction to Corollary \ref{cor:primero}.
Therefore $x_{\beta_i} \notin W_{\beta_j}$. Then the second statement of the lemma follows from Remark \ref{Rem:Walpha}.
\ep

We prove now the main result of this section.

\begin{theorem}\label{Thm:convexorder}
Keep the notation above. The order $\beta_1 < \beta_2 < \cdots < \beta_M$ on $\de^V_+$ is convex.
\end{theorem}
\bp
Each $W_{\beta_i}$ corresponds with one $F(w_i)$. As we have a chain as in previous lemma, by Theorem \ref{Thm:HSadaptado}
we have $w_1 \geq_D w_2 \geq_D \ldots \geq_D w_M$.

As the $w_i$'s are pairwise different, we have a chain of maximal length, and by Remark \ref{Rem:maximallength} there exists
a reduced expression of the longest element $\omega^V=\id_V s_{i_M} \cdots s_{i_1} $ such that $w_k=\id_V s_{i_M} \cdots
s_{i_k}$ for each $1 \leq k \leq M$.

We will prove by induction (descending on $j$) that $\beta_j= s_{i_M} \cdots s_{i_{j+1}}(\alpha_j)$. If so, we conclude the
proof because of Theorem \ref{Theorem:equivalenciasconvexo}. For the first step, notice that $\cH_{w_M}=
\bq_{N_{\alpha_\theta}}(x_\theta)$ by Theorem \ref{Thm:HSadaptado}, and by Remark \ref{Rem:Walpha} we have $i_m=\theta$,
so we have the initial step.

Assume now that $k < M$ and $\beta_j= s_{i_M} \cdots s_{i_{j+1}}(\alpha_j)$ for $j=k+1, \cdots, M$. Call $\gamma= s_{i_M} \cdots
s_{i_{k+1}}(\alpha_k)$, so by inductive hypothesis we have
$$  \cH_{W_{\beta_{k+1}}} = \prod_{j=k+1}^M \bq_{N_{\beta_j}}(X^{\beta_j}) , \quad \cH_{W_{\beta_k}} = \bq_{N_{\gamma}}
(X^{\gamma}) \left( \prod_{j=k+1}^M \bq_{N_{\beta_j}}(X^{\beta_j}) \right).  $$
On the other hand, $\{ x_{\beta_M}^{n_M} \cdots x_{\beta_k}^{n_k}: \ 0 \leq n_j < N_{\beta_j} \}$ is a linearly independent
set of $W_{\beta_k}$, so
$$ \cH_{W_{\beta_k}} \geq \prod_{j=k}^M \bq_{N_{\beta_j}}(X^{\beta_j}) , $$
where the inequality between the series means that the inequality holds for all the corresponding coefficients. Looking at
the coefficient of $X^{\beta_k}$ we obtain that there exists an expression
$$\beta_k= n \gamma + \sum_{j=k+1}^M n_j \beta_j, \qquad n \in \N, n_j \in \N_0.$$
Note that $R_{w_k}=\Lambda^{V}_+=\{ \gamma, \beta_{k+1}, \ldots \beta_M \}$, so if we apply $w_k$ to the last equality,
we obtain that $w_l^{-1}(\beta_k) \in \de^{r_l \cdots r_M(V)}_m$. Therefore $\beta_k \in R_{w_k}$, and as $\beta_k \neq
\beta_j$ for all $j > k$, we conclude $\beta_k=\gamma$.
\ep

The next result is analogous to the one for the positive part of quantized enveloping algebras $U_q(\mathfrak{g})$ given
in \cite{Le}, and gives an inductive way to obtain the words $l_{\beta}$ for $\beta \in \de^V_+$.

\begin{cor}\label{Coro:caracterizacionLyndon}
For each $\beta \in \de^V_+$, $\beta \neq \alpha_1, \ldots, \alpha_\theta$,
\begin{equation}\label{caracterizacionpalabrasLyndon}
    l_{\beta}= \max \{ l_{\delta_1}l_{\delta_2}: \ \ \delta_1, \delta_2 \in \de^V_+, \ \delta_1+\delta_2=\beta, \ l_{\delta_1}<l_{\delta_2} \}.
\end{equation}
\end{cor}
\bp Any factor of an element of $G_{I(V)}$ is in $G_{I(V)}$ (see Subsection \ref{subsection:pbw}). If $l_\beta=uv$ is the
Shirshov decomposition of $l_\beta$, then there exist $\gamma_1, \gamma_2 \in \de^V_+$ such that $u=l_{\gamma_1} <
v=l_{\gamma_2}$ and $\beta= \gamma_1+\gamma_2$.

On the other hand, let $\delta_1, \delta_2 \in \de^V_+$ be such that $\delta_1+\delta_2=\beta$ and $l_{\delta_1}<
l_{\delta_2}$. By the previous theorem, $l_{\delta_1}<l_\beta <l_{\delta_2}$. If $l_\beta$ does not begin with
$l_{\delta_1}$, then $l_{\delta_1}u <l_\beta$ for every word $u$, so in particular $l_{\delta_1}l_{\delta_2}<l_\beta$.
If $l_\beta$ begins with $l_{\delta_1}$, then $l_\beta= l_{\delta_1}u$, where $u$ has degree $\delta_2$. Let $u=l_p l_{p-1}
\cdots l_1$ be its Lyndon decomposition. Therefore each $l_i \in G_{I(V)}$, so $u=l_{\beta_M}^{n_M} \cdots
l_{\beta_1}^{n_1}$ for some $n_i \in \N_0$. Let $k =\max\{j: n_j \neq 0\}$. As the order is strongly convex, $x_{\beta_k}
\geq x_{\delta_2}$; i.e. $l_{\beta_k} \geq l_{\delta_2}$, so $u \geq l_{\delta_2}$ and then $l_\beta=l_{\delta_1}u \geq
l_{\delta_1}l_{\delta_2}$. In any case, $l_\beta=l_{\delta_1}u \geq l_{\delta_1}l_{\delta_2}$.
\ep

Another consequence shows that the family of coideal subalgebras $W_\beta$ (which are in particular left $\bB(V)$-comodules)
behaves as a kind of modules of highest weight.

\begin{theorem}\label{Thm:baseWbeta} The set $B_k =\{ x_{\beta_M}^{n_M} \cdots x_{\beta_k}^{n_k}: \ 0 \leq n_j < N_{\beta_j} \}$ is a basis of $W_{\beta_k}$. Moreover, if $W_{\beta_k} = \oplus_{\alpha \in \N^\theta_0} W_{\beta_k}(\alpha)$ denotes the decomposition in the $\N^\theta_0$
homogeneous components, then $\dim W_{\beta_k}(\beta_k)=1$.
\end{theorem}
\bp
The first statement follows because $B_k$ is included in $W_{\beta_k}$, it is linearly independent and the Hilbert series
of the $\kk$-linear subspace spanned by $B_k$ coincides with the Hilbert series of $W_{\beta_k}$.

For the second statement, if $\sum_{i=1}^M n_i \beta_i= \beta_k$ for some $n_i \in \N_0$, then $n_i= \delta_{i,k}$ or there
exists $i <k$ such that $n_i >0$, by Theorem \ref{Theorem:equivalenciasconvexo}.
\ep

The first consequence of the description of coideal subalgebras $W_\alpha$ as in previous theorem is a new expression of the coproduct of hyperwords which we will use in next section. We set
\begin{align}\label{subconjunto de la base PBW}
    C_k &:=\{ x_{\beta_k}^{n_k} x_{\beta_{k-1}}^{n_{k-1}} \cdots x_{\beta_1}^{n_1}: \ 0 \leq n_j < N_{\beta_j} \},
    \\ D_k &:=\{ x_{\beta_M}^{n_M} \cdots x_{\beta_1}^{n_1}: \ 0 \leq n_j < N_{\beta_j}, \ \exists j \geq k \mbox{ such that }n_j \neq 0 \}.
\end{align}

\begin{lema}\label{Lemma:productoBiBk}
Let $a \in B_k-\{1\}$, $b \in B_l$, $l  \leq k$. Then $ab=0$ or $ab$ is spanned by elements of $B_l \cap D_k$.
\end{lema}
\bp If $l=k$, it follows directly. Assume then $l < k$ and write $b=b_1b_2$ with $b_1 \in B_k$ and $b_2 \in C_{k-1} \cap B_l$ (possibly $b_1=1$). Then $ab_1 \in W_{\beta_k}$, because $W_{\beta_k}$ is a subalgebra, so it is spanned by $B_k$. To end, just note that if $c \in B_k$, then $cb_2 \in B_l \cap D_k$.
\ep

We set also $ht(u):=\sum n_i$, if $u=x_{\beta_M}^{n_M} x_{\beta_{k-1}}^{n_{k-1}} \cdots x_{\beta_1}^{n_1}$.

\begin{lema}\label{Lemma:coproductoforma2}
Let $u = x_{\beta_k}^{n_k} \cdots x_{\beta_l}^{n_l} \in B_l-D_{k+1}$, $l \leq k$, be such that $n_k, n_l \neq 0$. Then,
$$ \de(u) \in \left( \bigoplus_{v \in B, \ w \in D_k \cap B_l} \kk \ v \ot w \right) \bigoplus  \left( \bigoplus_{v \in D_k, \ w \in B_l-D_{k}} \kk \ v \ot w \right). $$
\end{lema}
\bp We prove it by induction on the height. If $ht(u)=1$, $u=x_{\beta_i}$ for some $i$. Then, $\de(u) \in u \ot 1+ 1 \ot u + \bB(V) \ot W_{\beta_i}$, so the result follows.

Assume it holds for $ht(w) < n$, and $u=x_{\beta_k}^{n_k} \cdots x_{\beta_l}^{n_l}$ is such that $ht(u)=n$. Write $u=x_{\beta_k} w$, so by inductive hypothesis,
$$ \de(u) \in \left( \bigoplus_{v \in B, \ w \in D_s \cap B_l} \kk \ v \ot w \right) \bigoplus  \left( \bigoplus_{v \in D_s, \ w \in B_l-D_{s}} \kk \ v \ot w \right), $$
where $s=k-1$ if $n_k=1$, or $s=k$ if $n_k>1$. We calculate $\de(u)=\de(x_{\beta_k})\de(w)$. Using that the braiding is diagonal and Lemma \ref{Lemma:productoBiBk} we conclude that
$$ ( \de(x_{\beta_k})- x_{\beta_k} \ot 1 ) \de(w) \in \bigoplus_{v \in B, \ w \in D_k \cap B_l} \kk \ v \ot w. $$
Also, for any $v \in B$ we have $x_{\beta_k} v \in D_k$, because if $v \in B_k$ then $x_{\beta_k} v \in W_{\beta_k}$ and if $v \in B_i$ for $i < k$ then we apply Lemma \ref{Lemma:productoBiBk} again, and we conclude the proof.
\ep

\bigskip

\section{Presentation by generators and relations of Nichols algebras of diagonal type}\label{section:relations} \

In this section we use the convex order of a PBW basis of hyperletters to prove that, when the diagonal braiding is
symmetric, such PBW basis is orthogonal with respect to the bilinear form of Proposition \ref{Prop:formabilineal}.
This fact gives a way to obtain relations which holds in the Nichols algebras, even when the braiding is not symmetric.
We obtain then a presentation by generators and relations for any Nichols algebras of diagonal type whose root system is
finite considering a suitable set of relations.

\subsection{A general presentation}\label{subsection:presentation}

We continue with the setting fixed in Subsection \ref{subsection:coidealsubalgebras}. To begin with, we prove the
orthogonality of the PBW basis with respect to the bilinear form in Proposition \ref{Prop:formabilineal}. This result
extends \cite[Prop. 5.1]{An}, and the proof is very similar; anyway we rewrite it in this general setting.

\begin{prop}\label{Prop:basePBWortogonal}
Consider a PBW basis of $\bB(V)$ as above given by Lyndon words, and assume that the braiding matrix is symmetric.
Then the PBW basis is orthogonal with respect to the bilinear form in Proposition \ref{Prop:formabilineal}.
\end{prop}

\bdem We prove by induction on $k=\ell (u)+ \ell (v)$
that $(u | v)=0$, where $u \neq v$ are ordered products of PBW generators. If
$k=1$, then $u=1$, $v=x_j$ or $u=x_i$, $v=1$, for some $i, j \in \unon$, and
$(1|x_j)=(x_i|1)=0$.

Suppose the statement is valid when $\ell (u)+ \ell (v) <k$, and let $u \neq v$ be hyperwords such that
$\ell (u)+ \ell (v) = k$. If both are hyperletters, they have
different degrees $\alpha \neq  \beta \in \zt$, so $u=x_{\alpha}$,
$v=x_{\beta}$, and $(x_{\alpha} | x_{\beta})=0$, since the
homogeneous components are orthogonal for $(\cdot | \cdot )$.

Suppose that $u=x_{\alpha}$ and $v=x_{\beta_k}^{h_k} x_{\beta_{k-1}}^{h_{k-1}} \ldots
x_{\beta_i}^{h_i}$, for some $1 \leq i \leq k \leq M$ (we consider $h_k, h_l \neq 0$).
If they have different $\zt$-degree, they are orthogonal. Then, we
assume that $\alpha= \sum_{j=i}^{k}h_j \beta_j$, so $\beta_i < \alpha$ because the ordered root system is strongly convex
by Theorem \ref{Thm:convexorder}. Using Lemma
\ref{Lemma:coproductPBWelements} and \eqref{Prop:formabilineal},
\begin{align*}
(u | v) =& ( x_{\alpha} | w) (1 | x_{\beta_i}) + (1 |w)
(x_{\alpha} | x_{\beta_i}) \\ &+ \sum_{ l_1\geq \dots \geq l_p
>l_{\alpha}, \ l_i \in
        L}  ( x_{l_1,\dots ,l_p} | w) ( [l_1]_c \cdots [l_p]_c | x_{\beta_i} )
\end{align*}
where $v=wx_{\beta_i}$. Note that $(1 | x_{\beta_i})= (1 |w) =0$.
Also, $[l_1]_c \cdots [l_p]_c$ is a linear combination of greater
hyperwords of the same degree and an element of $I(V)$, so by
inductive hypothesis and the fact that $I(V)$ is the radical of
the bilinear form, we conclude $( [l_1]_c \cdots [l_p]_c | x_{\beta_i} ) =0$. Therefore $(u | v)=0$.

For the final case, we consider
$$u=x_{\beta_k}^{h_k}  \ldots
x_{\beta_i}^{h_i}, \ 1 \leq i \leq k \leq M, \quad v=x_{\beta_q}^{f_q} \ldots
x_{\beta_p}^{f_p}, \ 1 \leq p \leq q \leq M. $$
The bilinear form is symmetric, so we can assume $i \leq p$. By Lemma \ref{Lemma:coproductPBWelements} and
\eqref{bilinearprop2},
\begin{align*}
    (u | v) &= ( w | 1) (x_{\beta_i} | v) + \sum ^{f_p}_{j=0} \binom{ f_p }{ j } _{q_{\beta_p}} (w | x_{\beta_q}^{f_q} \ldots x_{\beta_{p-1}}^{f_{p-1}}x_{\beta_p}^{j} ) (x_{\beta_i} | x_{\beta_p}^{f_p-j})
        \\ & \quad + \sum_{ \substack{ l_1\geq \dots  \geq l_t >l_{\beta_p}, \ l_s \in L \\ 0\leq j \leq f_p } } (w | x_{l_1,\dots ,l_t}^{(j)}) (x_{\beta_i} | \left[l_1\right]_c \dots      \left[l_t\right]_c x_{\beta_p}^j )
\end{align*}
where $u=w x_{\beta_i}$. Note
that $( w | 1)=0$, and $\left[l_1\right]_c \dots \left[l_p\right]_c
x_{\beta_p}^j$ is a combination of hyperwords of the PBW basis greater or equal than it and an element of $I(V)$. Using induction hypothesis and the fact that $I(V)$ is the radical of this bilinear form, we conclude that $(x_{\beta_i} | \left[l_1\right]_c \dots      \left[l_p\right]_c x_{\beta_p}^j )=0$. As also $x_{\beta_i}$, $x_{\beta_p}^{f_p-j} $ are different elements of the PBW basis for
$f_p-j \neq 1$, we have that
\begin{equation}\label{formabilinealhiperpalabras}
(u | v)= (f_p )_{q_{\beta_p}} (w | x_{\beta_q}^{f_q} \ldots x_{\beta_{p-1}}^{f_{p-1}} x_{\beta_p}^{f_p-1} ) (x_{\beta_i} | x_{\beta_p}) .
\end{equation}
Then it is zero if $i < p$, but also if $i=p$, because in that case $w \neq x_{\beta_q}^{f_q} \ldots x_{\beta_{p-1}}^{f_{p-1}} x_{\beta_p}^{f_p-1}$ and we use induction hypothesis. \edem

\begin{cor}\label{Coro:normageneradores}
If $u= x_{\beta_M}^{n_M} \cdots x_{\beta_1}^{n_1}$, where $0 \leq n_j < N_{\beta_j}$, then
\begin{equation}\label{normaPBWgenerador}
    c_u:= (u|u)= \prod_{j=1}^M n_j!_{q_{\beta_j}} c_{x_{\beta_j}}^{n_j} \neq 0.
\end{equation}
\end{cor}
\bp
We check the equality by induction on $ht(w)$. If $ht(w)=1$, $w$ is an hyperletter. If we assume it holds for $ht(w)<k$, and
$ht(u)=k$, we use the orthogonality of the PBW basis and a calculation as \eqref{formabilinealhiperpalabras} for $v=u$ to
deduce \eqref{normaPBWgenerador} from the inductive hypothesis.

Such scalar is not zero because $u \notin I(V)$ and the PBW basis generates a $\kk$-linear complement to $I(V)$, the radical
of this bilinear form.
\ep

\begin{obs}\label{Rem:calculoformabilineal} Note that:
$$ (x_{\beta_i}x_{\beta_j} | u)= (x_{\beta_i} | u_{(1)})(x_{\beta_j} | u_{(2)}) =  d_{i,j} c_{x_{\beta_i}}c_{x_{\beta_j}}, $$
where $d_{i,j}$ is the coefficient of $x_{\beta_i} \ot x_{\beta_j}$ in the expression of $\de(u)$ as a linear combination of
elements of the PBW basis in both sides of the tensor product.
\end{obs}

We return to the general case where the braiding matrix is not necessarily symmetric. We obtain some relations and
prove then the presentation of Nichols algebras by generators and relations. To obtain these
relations is the key step to find the presentation in Theorem \ref{Thm:presentacion}. Note that $B_i \cap C_j$ is
the set of monotone hyperwords whose hyperletters are between $x_{\beta_i}$ and $x_{\beta_j}$,see Theorem 
\ref{Thm:baseWbeta} and the definition of $C_j$ in Subsection \ref{subsection:coidealsubalgebras}.
\medskip

Let $(W,d)$ be a braided vector space of diagonal type, $ \hat x_1, \ldots, \hat x_\theta$ a basis of $W$ and
$\hat q_{ij} \in \kk^\times$ such that $d(\hat x_i \otimes \hat x_j) = \hat q_{ij} \hat x_j \otimes \hat x_i$. Assume that
$\hat q_{ij}=\hat q_{ji}$ for all $1 \leq i,j \leq \theta$, and that $(V,c)$ and $(W,d)$ are twist equivalent:
$$ q_{ij}q_{ji}= \hat q_{ij} \hat q_{ji}, \quad q_{ii}= \hat q_{ii}, \qquad 1 \leq i \neq j \leq \theta.$$
We define $\hat x_\beta= [l_\beta]_d$: that is, the corresponding hyperletter to $l_\beta$, but where we change the braiding
$c$ by $d$. By Corollary \ref{Coro:caracterizacionLyndon} and the invariance of the root system under twist equivalence,
the set of all the $\hat x_\beta$, $\beta \in \Delta^V_+= \Delta^W_+$, is a set of generators of a PBW basis as in
Kharchenko's Theorem. If $u= x_{\beta_M}^{n_M} \cdots x_{\beta_1}^{n_1}$, then we denote
$\hat u= \hat x_{\beta_M}^{n_M} \cdots \hat x_{\beta_1}^{n_1}$.

Let $\sigma: \zt \times \zt \rightarrow \kk^\times$ the bilinear form given by
\begin{equation}\label{cociclo}
 \sigma(g_i,g_j) = \left\{ \begin{array}{lr} \hat q_{ij}q_{ij}^{-1}, & i \leq j \\ 1, & i>j \end{array} \right.
\end{equation}

By \cite[Prop. 3.9, Rem. 3.10]{AS1} there exists a linear isomorphism $\Psi: \bB(W) \rightarrow \bB(V)$ such that
$\Psi(\hat x_i)=x_i$ and for any $x \in \bB(W)_\alpha$, $y \in \bB(W)_\beta$, $\alpha, \beta \in  \N^\theta_0$,
\begin{align}
\Psi(xy) &= \sigma( \alpha, \beta) \Psi (x) \Psi(y), \label{propiedadPsi1} \\
\Psi( [x,y]_d ) &= \sigma( \alpha, \beta) [\Psi (x), \Psi(y) ]_d. \label{propiedadPsi2}
\end{align}

Define $t_{\alpha_i}=1$ for all $1 \leq i \leq \theta$, and inductively
$$ t_\beta= \sigma(\beta_1, \beta_2) t_{\beta_1} t_{\beta_2}, \qquad \Sh(l_\beta)= (l_{\beta_1}, l_{\beta_2}) . $$
Also for any $u= x_{\beta_M}^{n_M} \cdots x_{\beta_1}^{n_1}$ define
\begin{equation}\label{formulacoeficientes}
 f(u):= \prod_{1 \leq i<j \leq M} \sigma(\beta_j, \beta_i)^{n_in_j} \ \prod_{1 \leq i \leq M} \sigma(\beta_i,
\beta_i)^{\binom{n_i}{2}}t_{\beta_i}^{n_i}.
\end{equation}

\begin{lema}\label{Lemma:calculo coef}
For any $u= x_{\beta_M}^{n_M} \cdots x_{\beta_1}^{n_1}$, $\Psi(\hat u) = f(u) u$.
\end{lema}
\bp
We prove first by induction on $\ell(l_\beta)$, $\beta \in \Delta^V_+$, that $\Psi(\hat x_\beta)=t_\beta \ x_\beta$. It follows by definition when $\ell(l_\beta)=1$, i.e. when $\beta=\alpha_i$ for some $1\leq i \leq \theta$. Now assume it holds for $\ell(l_\gamma) < k$, and consider $\beta \in \Delta^V_+$ such that $\ell(l_\beta)=k$. Let $Sh(l_\beta)=(\beta_1, \beta_2)$. Then,
\begin{align*}
\Psi(\hat x_\beta) &= \Psi([\hat x_{\beta_1}, \hat x_{\beta_2}]_d) = \sigma(\beta_1, \beta_2) [\Psi (\hat x_{\beta_1}), \Psi ( \hat x_{\beta_2})]_c
\\ &= \sigma(\beta_1, \beta_2) t_{\beta_1} t_{\beta_2}[ x_{\beta_1}, x_{\beta_2}]_c = t_\beta \ x_\beta,
\end{align*}
by \eqref{propiedadPsi2} and inductive hypothesis.

Now we prove that $\Psi(\hat u) = f(u) u$ by induction on $ht(u)$. Note that if $ht(u)=1$ it reduces to $\Psi(\hat x_\beta)=t_\beta \ x_\beta$. Assume now that it holds for $ht(v)< N$, and consider $u= x_{\beta_M}^{n_M} \cdots x_{\beta_k}^{n_k}$ such that $ht(u)=N$ and $n_k > 0$. Call $v=  x_{\beta_M}^{n_M} \cdots x_{\beta_k}^{n_k-1}$. Then,
\begin{align*}
\Psi(\hat u) &= \sigma \left( (n_k-1)\beta_k +\sum_{i=k+1}^M n_i \beta_i, \beta_k \right)\Psi( \hat v) \Psi( \hat x_{\beta_k} )
\\ &= \left( \prod_{i=k+1}^M \sigma(\beta_i,\beta_k)^{n_i} \right) \sigma(\beta_k, \beta_k)^{n_k-1} f(v)v \ t_{\beta_k} x_{\beta_k}= f(u) u,
\end{align*}
by \eqref{propiedadPsi1} and inductive hypothesis.
\ep

We define for $ 1\leq i<j \leq \theta$ and $u= x_{\beta_M}^{n_M} \cdots x_{\beta_1}^{n_1}$,
\begin{equation}\label{escalares}
c_{i,j}^u := \frac{f(u) \ (\hat x_{\beta_i} \hat x_{\beta_j} | \hat u)}{\sigma(\beta_i, \beta_j) t_{\beta_i} t_{\beta_j} c_{\hat u}},
\end{equation}
where $(\cdot|\cdot)$ denotes the bilinear form associated to $(W,d)$, and $c_{\hat u}$ is computed as in Corollary \ref{Coro:normageneradores}. Note that if $(q_{ij})$ is symmetric and we consider $q_{ij}= \hat q_{ij}$, then $\sigma(\alpha, \beta)=1$ for all $\alpha, \beta \in \zt$ and then $f(u)=1$ for any $u$. In consequence, $c_{i,j}^u= (x_{\beta_i}x_{\beta_j} | u) c_u^{-1}$.

We obtain a first set of relations for our presentation.

\begin{lema}\label{Lemma:otrasrelaciones} Let $1 \leq i<j \leq M$ be such that $l_{\beta_i} l_{\beta_j} \neq l_{\beta_k}$ for all $k$, and $\Sh(l_{\beta_i} l_{\beta_j}) =( l_{\beta_i}, l_{\beta_j} )$, and $c_{i,j}^u \in \kk$ as above. Then,
\begin{equation}\label{quantumSerregeneralizadas}
    \left[ x_{\beta_i}, x_{\beta_j} \right]_c= \sum_{u \in B_i \cap C_j- \{ x_{\beta_j}x_{\beta_i} \}: \ \ \deg u= \beta_i+\beta_j} c_{i,j}^u \ u.
\end{equation}
\end{lema}
\bp Assume that the braiding is symmetric. As $l_{\beta_i} l_{\beta_j} \neq l_{\beta_k}$ for all $k$, and $\Sh(l_{\beta_i} l_{\beta_j}) =( l_{\beta_i}, l_{\beta_j} )$, $[l_{\beta_i} l_{\beta_j}]_c = [x_{\beta_i}, x_{\beta_j}]_c = x_{\beta_i} x_{\beta_j} - \chi (\beta_i, \beta_j) x_{\beta_j} x_{\beta_i}$ is a linear combination of greater monotone hyperwords by Corollary \ref{cor:primero}.

As $x_{\beta_i}x_{\beta_j} \in W_{\beta_i}$, it is a linear combination of elements in $B_i$ by Theorem \ref{Thm:baseWbeta}.  Also, $\bB(V)$ is $\N^\theta_0$-graded, so this linear combination is over elements of $B_i$ of the degree $\beta_i+\beta_j$. Moreover, if $c_{i,j}^u \neq 0$ for $u = x_{\beta_k}^{n_k} \cdots x_{\beta_l}^{n_l}$, $l \leq k$, such that $n_k, n_l \neq 0$, then $x_{\beta_i} \ot x_{\beta_j}$ appears in the expression of $\de(u)$ by Remark \ref{Rem:calculoformabilineal}. Note that $ x_{\beta_i} \ot x_{\beta_j} \notin D_k \ot (B_l-D_{k})$, because $i <j$. By Lemma \ref{Lemma:coproductoforma2}, we have $x_{\beta_j} \in B_k$, so $j \geq k$, and $u \in C_j$.

The explicit formula of the coefficients comes from Proposition \ref{Prop:basePBWortogonal}.

If we want to compute $c_{i,j}^{x_{\beta_j}x_{\beta_i} }$, we have to calculate the coefficient of $x_{\beta_i} \ot x_{\beta_j}$ in $\de(x_{\beta_j}x_{\beta_i})$, because of Remark \ref{Rem:calculoformabilineal} and the formula $c_{x_{\alpha_j}x_{\alpha_i}} = c_{x_{\alpha_i}} c_{x_{\alpha_j}}$. This coefficient is $\chi(\beta_j,\beta_i)$, but as the braiding matrix is symmetric, $\chi(\beta_j,\beta_i)=\chi(\beta_i,\beta_j)$. Therefore we conclude the proof when the matrix braiding is symmetric.

When the braiding is not symmetric, we use the linear isomorphism $\Psi$ considered previously to reduce the computation to
the symmetric case. Then,
\begin{align*}
0 &= \Psi \left( \left[ \hat x_{\beta_i}, \hat x_{\beta_j} \right]_d - \sum  (\hat x_{\beta_i}\hat x_{\beta_j} | \hat u) c_{\hat u} ^{-1} \ \hat u \right)
\\ &= \sigma(\beta_i, \beta_j) t_{\beta_i} t_{\beta_j} [ x_{\beta_i}, x_{\beta_j} ]_c - \sum  (\hat x_{\beta_i}\hat x_{\beta_j} | \hat u) c_{\hat u} ^{-1} f(u) u,
\end{align*}
by \eqref{propiedadPsi2} and Lemma \ref{Lemma:calculo coef}, so \eqref{quantumSerregeneralizadas} holds in $\bB(V)$.
\ep

\begin{cor}\label{Coro:otrasrelaciones}
Assume that $i,j$ are as in Lemma \ref{Lemma:otrasrelaciones}, and $\beta_i+\beta_j= \sum_{k=i}^j n_k \beta_k$, $n_k \in \N_0$ if and only if $n_i=n_j=1$, $n_k=0$ for $k \neq i,j$. Then,
\begin{equation}\label{quantumSerregeneralizadas2}
[x_{\beta_i}, x_{\beta_j} ]_c=0.
\end{equation}
\end{cor}
\bp It follows from the previous proposition.
\ep

Now we extend \cite[Cor. 5.2]{An}. Recall that $N_\beta= \ord(q_\beta)= h(x_\beta)$.

\begin{lema}\label{Lemma:heigthgenerators}
If $\beta \in \de^V_+$ and $N_\beta$ is finite, then
\begin{equation}\label{powerrootvector}
x_{\beta}^{N_{\beta}}=0, \qquad \mbox{in }\bB(V).
\end{equation}
\end{lema}
\bdem Assume first that $(q_{ij})$ is symmetric. Consider $w= \widetilde{w} x_\beta^m$, where $\beta \in \Delta^+$ and $\widetilde{w}$ is a non-increasing product of hyperletters $x_\gamma$, $\gamma \in \Delta^+,  \gamma > \beta$ or $\widetilde{w}=1$. If $\beta >\alpha$,
\begin{align*}
    (x_{\alpha}^{N_{\alpha}} | w) &= ( x_{\alpha}^{N_{\alpha}-1} | 1) (x_{\alpha} | w) + \sum ^{m}_{i=0} \binom{ m }{ i } _{q_{\beta}} ( x_{\alpha}^{N_{\alpha}-1} | \widetilde{w} x_{\beta}^{i} ) (x_{\alpha} | x_{\beta}^{m-i})
        \\ & \quad + \sum_{ \substack{ l_1\geq \dots  \geq l_p >x_{\beta}, 0\leq j \leq m } } (x_{\alpha}^{N_{\alpha}-1} | x_{l_1,\dots ,l_p}^{(j)}) (x_{\alpha} | \left[l_1\right]_c \dots      \left[l_p\right]_c x_{\beta}^j )=0,
\end{align*}
where we use that $( x_{\alpha}^{N_{\alpha}-1} | 1)= (x_{\alpha} | x_{\beta}^{m-i})=(x_{\alpha} | \left[l_1\right]_c \dots      \left[l_p\right]_c x_{\beta}^j )=0$ by the orthogonality of the PBW basis.

If $\beta \leq \alpha$, then
\begin{align*}
    (x_{\alpha}^{N_{\alpha}} | w) &= ( 1 | \widetilde{w} x_{\beta}^{m-1}) (x_{\alpha}^{N_{\alpha}} | x_{\beta}) + \sum ^{N_{\alpha}}_{i=1} \binom{ N_{\alpha} }{ i } _{q_{\alpha}} ( x_{\alpha}^i | \widetilde{w} x_{\beta}^{m-1} ) ( x_{\alpha}^{N_{\alpha}-i} | x_{\beta})
        \\ & \quad + \sum_{ \substack{ l_1\geq \dots  \geq l_p >x_{\alpha}, 0\leq j \leq N_{\alpha} } } ( x_{l_1,\dots ,l_p}^{(j)}  | \widetilde{w} x_{\beta}^{m-1} ) ( \left[l_1\right]_c \dots      \left[l_p\right]_c x_{\alpha}^j | x_{\beta} )
\end{align*}
where we use that $q_{\alpha} \in \G_{N_{\alpha}}$, the orthogonality of the PBW basis and the fact that $N\beta \notin \Delta^+$ if $N>1$ (so $(x_{\alpha}^{N_{\alpha}} | x_{\beta})=0$).

Therefore $(x_{\alpha}^{N_{\alpha}} |v)=0$ for all $v$ in the PBW basis. Also $(I(V) | x_{\alpha}^{N_{\alpha}})=0$, because it is the radical of this bilinear form, so $(T(V) | x_{\alpha}^{N_{\alpha}})=0$, and then $x_{\alpha}^{N_{\alpha}} \in I(V)$. That is, we have
$x_{\alpha}^{N_{\alpha}}=0$  in $\bB(V)$.
\medskip

For the general case, we recall that a diagonal braiding is twist
equivalent to a braiding with a symmetric matrix, see
\cite[Theorem 4.5]{AS1}. Also, there exists a linear isomorphism
between the corresponding Nichols algebras. The corresponding
$x_{\alpha}$ are related by a non-zero scalar, because they are
an iteration of braided commutators between the hyperwords.
\edem

Before proving the main result of this section, we need another technical lemma.

\begin{lema}\label{Lemma:productopalabras}
Let $\bB$ be a quotient of $T(V)$ such that relations \eqref{quantumSerregeneralizadas} hold. Then for any $i < j$, $x_{\beta_i} x_{\beta_j}$ can be written as a linear combination of monotone hyperwords greater than $x_{\beta_i}$, whose hyperletters are $x_{\beta_k}$, $i \leq k\leq j$.
\end{lema}
\bp
It is similar to the proof of Theorem \ref{Theo:corcheteentreLyndon}, see \cite[Thm. 10]{R2}. Set for each $n \geq 2$,
$$ L_n:= \{ (x_{\beta_i}, x_{\beta_j}) : \ i< j, \ \ell(l_{\beta_i})+ \ell(l_{\beta_j})=n \}. $$
We order $L_k$ as follows: $(x_{\beta_i}, x_{\beta_j}) < (x_{\beta_k}, x_{\beta_m})$ if $l_{\beta_i}l_{\beta_j}) < l_{\beta_k}l_{\beta_m}$, or $l_{\beta_i}l_{\beta_j}=l_{\beta_k}l_{\beta_m}$ and $l_{\beta_i}< l_{\beta_k}$.

We prove the statement by induction on $n=\ell(x_{\beta_i}) + \ell(x_{\beta_j})$, and then by induction on the previous order on $L_n$. When $n=2$, then $\beta_i, \beta_j$ are simple, and $[x_i, x_j]_c=x_{\alpha_i+\alpha_j}$ or $[x_i, x_j]_c=0$ in $\bB$.

Fix then a pair $ (x_{\beta_i}, x_{\beta_j}) \in L_n$ and assume that the statement holds for $(x_{\beta_k}, x_{\beta_m}) \in L_n$, $(x_{\beta_i}, x_{\beta_j}) > (x_{\beta_k}, x_{\beta_m})$, and for $ (x_{\beta_k}, x_{\beta_m}) \in L_{n'}$, $n' <n$. If $\Sh(l_{\beta_i}l_{\beta_j})=(l_{\beta_i},l_{\beta_j})$ then the assertion holds because
\begin{itemize}
  \item if $l_{\beta_i}l_{\beta_j}= l_{\beta_k}$ for some $k$, necessarily (by the definition of the order) $i< k <j$ and $[x_{\beta_i}, x_{\beta_j}]_c= x_{\beta_k}$,
  \item if not, it holds because we assume \eqref{quantumSerregeneralizadas}.
\end{itemize}
If $\Sh(l_{\beta_i}l_{\beta_j}) \neq (l_{\beta_i},l_{\beta_j})$, let $\Sh(l_{\beta_i})=(l_{\beta_p}, l_{\beta_q})$, so $x_{\beta_i}= [x_{\beta_p}, x_{\beta_q}]_c$. Therefore $l_{\beta_q} < l_{\beta_j}$ (see Subsection \ref{subsection:pbw}). By \eqref{idjac},
$$ [x_{\beta_i}, x_{\beta_j}]_c = [x_{\beta_p}, [x_{\beta_q}, x_{\beta_j}]_c ]_c -\chi(\beta_p,\beta_q) x_{\beta_q} [x_{\beta_p}, x_{\beta_j}]_c+ \chi(\beta_q, \beta_j) [x_{\beta_p}, x_{\beta_j}]_c x_{\beta_q}. $$
We apply induction hypothesis and express $[x_{\beta_q}, x_{\beta_j}]_c$ as a linear combination of monotone hyperwords whose hyperletters are between $x_{\beta_q}$ and $x_{\beta_j}$. By \eqref{der} and inductive hypothesis, we express $[x_{\beta_p}, [x_{\beta_q}, x_{\beta_j}]_c ]_c$ as a linear combination of monotone hyperwords whose letters are between $x_{\beta_i}$ and $x_{\beta_j}$. It is important here the order in $L_n$, because in such linear combination can appear a single hyperletter $x_{\beta_k}$, which by hypothesis is between $x_{\beta_q}$ and $x_{\beta_j}$, and so $(l_{\beta_i}, l_{\beta_j}) > (l_{\beta_p}, l_{\beta_k})$.

We use also inductive hypothesis to express $[x_{\beta_p}, x_{\beta_j}]_c$ as a linear combination of hyperwords whose hyperletters are between $x_{\beta_p}$ and $x_{\beta_j}$. As in the previous step we can reorder the hyperletters in order to find the desired expression by inductive hypothesis.
\ep

Now we are ready to prove the main result of this work.

\begin{theorem}\label{Thm:presentacion}
Let $(V,c)$ be a finite-dimensional braided vector space of diagonal type such that $\de^V_+$ is finite. Let $x_1, \cdots, x_\theta$ be a basis of $V$ such that $c(x_i \ot x_j)= q_{ij} x_j \ot x_i$, where $(q_{ij}) \in (\kk^\times)^{\theta \times \theta}$ is the braiding matrix, and let $\{ x_{\beta_k} \}_{\beta_k \in \de^V_+}$ be the associated set of hyperletters.

Then $\bB(V)$ is presented by generators $x_1, \ldots, x_\theta$, and relations
\begin{align}
x_{\beta}^{N_{\beta}}& =0, \qquad \beta \in \Delta^V_+, \ \ord(q_\beta)=N_\beta < \infty, \label{powerrootvector1}
\\ \left[ x_{\beta_i}, x_{\beta_j} \right]_c &= \sum_{u \in B_i \cap C_j- \{ x_{\beta_j}x_{\beta_i} \}: \ \ \deg u= \beta_i+\beta_j} c_{i,j}^u \ u, \label{quantumSerregeneralizadas1}
\\ & 1 \leq i<j \leq M, \ \Sh(l_{\beta_i} l_{\beta_j}) =( l_{\beta_i}, l_{\beta_j} ), \ l_{\beta_i} l_{\beta_j} \neq l_{\beta_k}, \forall k, \nonumber
\end{align}
where $c_{i,j}^u$ are as in \eqref{escalares}. Moreover, $\{ x_{\beta_M}^{n_M} \cdots x_{\beta_1}^{n_1}: \ 0 \leq n_j < N_{\beta_j} \}$ is a basis of $\bB(V)$.
\end{theorem}
\bp
The statement about the basis follows by Kharchenko's theory on PBW bases (see Subsection \ref{subsection:pbw}) and the definition of $\de^V_+$ (see Subsection \ref{subsection:weylgroupoid}), where the hyperletters $x_{\beta_k}$ are univocaly determined by Corollary \ref{Coro:caracterizacionLyndon}.

Let $\bB :=T(V)/I$, where $I$ is the ideal of $T(V)$ generated by \eqref{powerrootvector1}, \eqref{quantumSerregeneralizadas1}: by Lemmata \ref{Lemma:otrasrelaciones} and \ref{Lemma:heigthgenerators}, $I \subseteq I(V)$, so the projection $\pi: T(V) \twoheadrightarrow \bB(V)$ induces canonically a projection $\pi': \bB \twoheadrightarrow \bB(V)$. Let $W$ be the subspace of $\bB$ spanned by $B$, where $B$ is the PBW basis of $\bB(V)$; $1 \in W$. For each pair $1 \leq i \leq j \leq M$, we set $W_{i,j}$ the subspace of $W$ spanned by $B_i \cap C_j$.

We assert that
\begin{equation}\label{afirmacionWij}
x_{\beta_k} W_{i,j} \subset W_{\min \{i,k\}, \max \{j,k\} }.
\end{equation}
We shall prove it by induction on $k$. When $k=M$, fix $i \leq j$. For each $w \in B_i \cap C_j$, we have that $x_{\beta_M}w \in B_i \cap C_M = B_i$, or $x_{\beta_M}w=0$ if $j=M$, $N_M < \infty$ and $w$ begins with $x_{\beta_M}^{N_M-1}$, so $x_{\beta_M}W_{i,j} \subset W_{i,M}$.

Now assume that \eqref{afirmacionWij} holds for all $l > k$ and all $i \leq j$. We argue by induction on $j$. If $i \leq j \leq k$, for each $w \in B_i \cap C_j$, we have that $x_{\beta_k}w \in B_i \cap C_k$ or $x_{\beta_k}w=0$ as in the initial step, so $x_{\beta_k}W_{i,j} \subset W_{i,k}$. Now assume $j > k$, and consider $w \in B_i \cap C_j$; it is enough to prove that $ x_{\beta_k} w \in W_{\min \{i,k\}, j}$. Moreover, we can assume $w=x_{\beta_j}w'$ for some monotone hyperword $w'$ in $W_{i,j}$ (if $w$ begins with another hyperletter $x_{\beta_l}$, $l < j$, we consider $w \in W_{i,l} \subset W_{i,j}$). By Lemma \ref{Lemma:productopalabras}, we can write $x_{\beta_k}x_{\beta_j}$ as a linear combination of monotone hyperwords whose hyperletters belong to $B_k \cap C_j$. Therefore the result follows by the inductive hypothesis: any of these hyperwords has no letters $x_{\beta_k}$'s and we use the first inductive hypothesis (it holds for all $l >k$), or it ends with hyperletters $x_{\beta_k}$'s and we write $x_{\beta_k}w'$ as a linear combination of hyperwords in $B_{\min\{i,k\}} \cap C_j$ by the second inductive hypothesis.

In this way we prove that $W$ is a left ideal which contains 1, so $W= \bB$. But then the projection $\pi'$ is an isomorphism, and $\bB=\bB(V)$.
\ep

\begin{obs}\label{Rem:quantumSerre}
Recall that we have defined for $i, j\in \{ 1, \ldots, \theta \}$,
$$ m_{ij}:= \max \{m: (\ad_c x_i)^m x_j \neq 0 \},$$
see \eqref{mij}, and then $m\alpha_i+\alpha_j \in \de^V_+$ iff $0 \leq m \leq m_{ij}$. Moreover assume $i<j$. Then $x_{m\alpha_i+\alpha_j}= (\ad_c x_i)^m x_j$, and a pair as in Corollary \ref{Coro:otrasrelaciones} is $(x_i, x_i^{m_{ij}} x_j)$, so such corollary implies the well-known quantum Serre relation in $\bB(V)$: $(\ad_c x_i)^{m_{ij}+1} x_j=0$. If $i> j$, then the pair changes to $(x_j x_i^{m_{ij}}, x_i)$, but then $0=[x_{m\alpha_i+\alpha_j}, x_i]_c= a (\ad_c x_i)^{m_{ij}+1} x_j$ for some $a \in \kk^\times$. In any case we have $(\ad_c x_i)^{m_{ij}+1} x_j=0$.

This shows that the set of relations \eqref{quantumSerregeneralizadas}, \eqref{powerrootvector} is not minimal: if $\ord q_{ii}=m_{ij}+1$, then $x_i^{m_{ij}+1}$ is one of the relations \eqref{powerrootvector}, and then $(\ad_c x_i)^{m_{ij}+1} x_j$ belongs to the ideal generated by $x_i^{m_{ij}+1}$.
\end{obs}

\section{Explicit presentations by generators and relations of some Nichols algebras of diagonal type}\label{section:examples}

We shall apply the previous theory about how to obtain a PBW basis (Corollary \ref{Coro:caracterizacionLyndon}) and
a presentation of the corresponding Nichols algebra (Theorem \ref{Thm:presentacion}) in some concrete examples.

\subsection{Examples when $\dim V=3$}

We consider the Weyl equivalence classes 9, 10, 11 in \cite[Table 2]{H2}. We fix the following notation: let $q,r,s \in
\kk^\times$ be such that $qrs=1$. Set $M, N,P \in \N$ as the orders of these scalars, if they are finite. Such Weyl
equivalence class includes the following generalized Dynkin diagrams:
\begin{itemize}
 \item $\xymatrix{ \circ^{q} \ar@{-}[r]^{q^{-1}} & \circ^{-1} \ar@{-}[r]^{r^{-1}} & \circ^{r} }$,
 \item $\xymatrix{ \circ^{q} \ar@{-}[r]^{q^{-1}} & \circ^{-1} \ar@{-}[r]^{s^{-1}} & \circ^{s} }$,
 \item $\xymatrix{ \circ^{r} \ar@{-}[r]^{r^{-1}} & \circ^{-1} \ar@{-}[r]^{s^{-1}} & \circ^{s} }$,
 \item $\xymatrix{ & \circ^{-1} \ar@{-}[rd]^{r} & \\ \circ^{-1} \ar@{-}[ru]^{q} \ar@{-}[rr]^{s} &  & \circ^{-1} }$.
\end{itemize}
Notice that 10, 11 are particular cases of 9 when $q=r$, $q=r=s \in \G_3$, respectively. Also the second and the third
diagrams are analogous to the first one, so it is enough to obtain the presentation for the first and the last braidings.

If $i<j$, $l_{\alpha_i+\alpha_j}=x_ix_j$, so $x_{\alpha_i+\alpha_j}= [x_i, x_j]_c= (\ad_c x_i) x_j$. Also,
\begin{equation*}
l_{\alpha_1+\alpha_2+\alpha_3}= \left\{ \begin{array}{lr} x_1x_2x_3 \mbox{ if }(\ad_c x_1)x_3=0:
&  x_{\alpha_1+\alpha_2+\alpha_3}= [x_1, x_{\alpha_2+\alpha_3}]_c;
\\ x_1x_3x_2 \mbox{ if }(\ad_c x_1)x_3 \neq 0:  &  x_{\alpha_1+\alpha_2+\alpha_3}= [x_{\alpha_1+\alpha_3}, x_2 ]_c;
\end{array} \right.
\end{equation*}
When $(\ad_c x_1)x_3=0$, we have also
$$ l_{\alpha_1+2\alpha_2+\alpha_3}=x_1x_2x_3x_2 : \qquad x_{\alpha_1+2\alpha_2+\alpha_3}= [x_{\alpha_1+\alpha_2+\alpha_3},
x_2]_c $$

\begin{prop}\label{Prop:ejemplo1}
Let $(V,c)$ be a braided vector space such that $\dim V=3$, and the corresponding generalized Dynkin diagram is
$$\xymatrix{ \circ^{q} \ar@{-}[r]^{q^{-1}} & \circ^{-1} \ar@{-}[r]^{r^{-1}} & \circ^{r} }. $$
Then $\bB(V)$ is presented by generators $x_1,x_2,x_3$, and relations
\begin{align}
    & x_1^M=x_2^2=x_3^N=x_{\alpha_1+2\alpha_2+\alpha_3}^P=0, \label{power-ej1}
    \\ & (\ad_c x_1)^2x_2= (\ad_c x_3)^2 x_2 = (\ad_c x_1)x_3=0, \label{qSerre-ej1}
    \\ & [ x_{\alpha_1+\alpha_2}, x_{\alpha_1+\alpha_2+\alpha_3}]_c = [ x_{\alpha_1+\alpha_2+\alpha_3},  x_{\alpha_2+\alpha_3} ]_c =0. \label{genqSerre-ej1}
\end{align}
Moreover, $\bB(V)$ has a PBW basis as follows:
 \begin{multline*}
  \left\{ x_3 ^{n_3} \ x_{\alpha_2+\alpha_3}^{n_{23}} \ x_2^{n_2} \  x_{\alpha_1+2\alpha_2+\alpha_3}^{n_{1232}} \ x_{\alpha_1+\alpha_2+\alpha_3}^{n_{123}} \ x_{\alpha_1+\alpha_2}^{n_{12}} \ x_1^{n_1} : \right.
 \\  \left. 0 \leq n_1 < M, \ 0 \leq n_2 < N, \ 0 \leq n_{1232} < P, \  n_{12}, n_{123}, n_2, n_{23} \in \{ 0,1 \} \right\}.
 \end{multline*}
If $M, N, P < \infty$, then $\dim \bB(V)= 16 MNP$.
\end{prop}
\bp For this case,
$$ \de^V_+= \{ \alpha_3, \alpha_2+\alpha_3, \alpha_2, \alpha_1+2\alpha_2+\alpha_3, \alpha_1+\alpha_2+\alpha_3, \alpha_1+\alpha_2, \alpha_1 \}.  $$
Therefore we obtain $l_\beta$, $\beta \in \de^V_+$, easily from Corollary \ref{Coro:caracterizacionLyndon}.

By Remark \ref{Rem:quantumSerre}, we consider the relations $$ (\ad_c x_1)^2x_2= (\ad_c x_3)^2 x_2 = (\ad_c x_1)x_3=0, $$
because $(\ad_c x_2)^2x_1,(\ad_c x_2)^2x_3$ follows from $x_2^2=0$.

We have the following decompositions:
\begin{align*}
\Sh(l_{\alpha_1+\alpha_2} l _{\alpha_1+\alpha_2+\alpha_3})&= (l_{\alpha_1+\alpha_2}, l_{\alpha_1+\alpha_2+\alpha_3}), \\
\Sh(l_{\alpha_1+\alpha_2+\alpha_3} l _{\alpha_2+\alpha_3})&= (l_{\alpha_1+\alpha_2+\alpha_3}, l_{\alpha_2+\alpha_3}).
\end{align*}
Relations \eqref{genqSerre-ej1} then follows by Corollary \ref{Coro:otrasrelaciones}.

Also $\Sh(l_{\alpha_1} l _{\alpha_1+\alpha_2+\alpha_3})= (l_{\alpha_1}, l_{\alpha_1+\alpha_2+\alpha_3})$, so $$ [x_1, x_{\alpha_1+\alpha_2+\alpha_3}]_c=0. $$
Note that $x_{\alpha_1+\alpha_2+\alpha_3}= [x_{\alpha_1+\alpha_2}, x_3 ]_c$ by $(\ad_c x_1)x_3=0$ and the identity \eqref{idjac}. Therefore, this relation is redundant because of \eqref{idjac}, $x_1^2=0$. The same holds for relation $[x_{\alpha_1+\alpha_2+\alpha_3}, x_3]_c=0$, coming from the decomposition
$\Sh(l _{\alpha_1+\alpha_2+\alpha_3}l_{\alpha_3})= (l_{\alpha_1+\alpha_2+\alpha_3}, l_{\alpha_3})$.

Also $\Sh(l_{\alpha_1} l _{\alpha_1+2\alpha_2+\alpha_3})= (l_{\alpha_1}, l_{\alpha_1+2\alpha_2+\alpha_3})$, so by Lemma \ref{Lemma:otrasrelaciones} there exists $a \in \kk$ such that:
$$ [ x_1, x_{\alpha_1+2\alpha_2+\alpha_3}]_c = a \ x_{\alpha_1+\alpha_2+\alpha_3} x_{\alpha_1+\alpha_2}. $$
This relation is also redundant:
\begin{align*}
    [ x_1, x_{\alpha_1+2\alpha_2+\alpha_3}]_c = &  [[ x_1, x_{\alpha_1+\alpha_2+\alpha_3}]_c, x_2 ]_c +q_{11}q_{12}q_{13} x_{\alpha_1+\alpha_2+\alpha_3} x_{\alpha_1+\alpha_2}
\\ & \ \ - q_{12}q_{22}q_{32} x_{\alpha_1+\alpha_2}  x_{\alpha_1+\alpha_2+\alpha_3}
\\ = & q_{11}q_{12}q_{13}(1-s) x_{\alpha_1+\alpha_2+\alpha_3} x_{\alpha_1+\alpha_2}.
\end{align*}
where we use \eqref{idjac} and the previous relations.

We have finally
\begin{align*}
\Sh(l_{\alpha_1+\alpha_2} l _{\alpha_1+2\alpha_2+\alpha_3}) &= (l_{\alpha_1+\alpha_2}, l_{\alpha_1+2\alpha_2+\alpha_3}), \\ \Sh( l _{\alpha_1+\alpha_2+\alpha_3} l_{\alpha_1+2\alpha_2+\alpha_3} ) &= ( l _{\alpha_1+\alpha_2+\alpha_3} , l_{\alpha_1+2\alpha_2+\alpha_3}),
\end{align*}
which give place to the following relations:
$$ [ x_{\alpha_1+\alpha_2}, x_{\alpha_1+2\alpha_2+\alpha_3}]_c = [ x_{\alpha_1+\alpha_2+\alpha_3} , x_{\alpha_1+2\alpha_2+\alpha_3} ]_c=0 $$
These relations also follow by the previous ones using \eqref{idjac}.

We can prove in the same way that $$ x_{\alpha_1+\alpha_2}^2, \ x_{\alpha_2+\alpha_3}^2, \ x_{\alpha_1+\alpha_2+\alpha_3}^2 $$
are redundant relations too. The proposition follows then by Theorem \ref{Thm:presentacion}, where we omit some redundant relations.
\ep

\begin{prop}\label{Prop:ejemplo2}
Let $(V,c)$ be a braided vector space such that $\dim V=3$, and the corresponding generalized Dynkin diagram is
$$
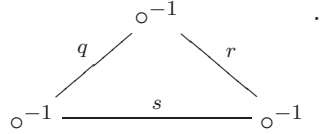
\xymatrix{ & \circ^{-1} \ar@{-}[rd]^{r} & \\ \circ^{-1} \ar@{-}[ru]^{q} \ar@{-}[rr]^{s} &  & \circ^{-1} }. $$
Then $\bB(V)$ is presented by generators $x_1,x_2,x_3$, and relations
\begin{align}
    & x_1^2=x_2^2=x_3^2=x_{\alpha_1+\alpha_2+\alpha_3}^2=0, \label{power-ej2}
    \\ & x_{\alpha_1+\alpha_2}^M= x_{\alpha_2+\alpha_3}^N = x_{\alpha_1+\alpha_3}^P=0, \label{powernot2-ej2}
    \\ & [ x_{\alpha_i+\alpha_j}, x_{\alpha_i+\alpha_k}]_c = 0, \qquad \{i,j,k \}= \{1,2,3 \}. \label{genqSerre-ej2}
    \\ & [ x_1, x_{\alpha_2+\alpha_3}]_c = \frac{1-s}{q_{23} (1-r)} x_{\alpha_1+\alpha_2+\alpha_3} + q_{12} (1-s)
x_2 x_{\alpha_1+\alpha_3}. \label{genqSerre2-ej2}
\end{align}
Moreover, $\bB(V)$ has a PBW basis as follows:
 \begin{multline*}
  \left\{ x_3 ^{n_3} \ x_{\alpha_2+\alpha_3}^{n_{23}} \ x_2^{n_2} \  x_{\alpha_1+\alpha_3}^{n_{13}} \ x_{\alpha_1+\alpha_2+\alpha_3}^{n_{123}} \ x_{\alpha_1+\alpha_2}^{n_{12}} \ x_1^{n_1} : \right.
 \\  \left. 0 \leq n_{12} < M, \ 0 \leq n_{23} < N, \ 0 \leq n_{13} < P, \  n_1, n_{123}, n_2, n_3 \in \{ 0,1 \} \right\}.
 \end{multline*}
If $M, N, P < \infty$, then $\dim \bB(V)= 16 MNP$.
\end{prop}
\bp Again we obtain $l_\beta$, $\beta \in \de^V_+$, easily from Corollary \ref{Coro:caracterizacionLyndon}, because
$$ \de^V_+= \{ \alpha_3, \alpha_2+\alpha_3, \alpha_2, \alpha_1+\alpha_3, \alpha_1+\alpha_2+\alpha_3, \alpha_1+\alpha_2, \alpha_1 \}.  $$

By Remark \ref{Rem:quantumSerre}, all the quantum Serre relations $(\ad_c x_i)^2 x_j=0$, $i \neq j$, follow from $x_i^2=0$, $i=1,2,3$.

We have the decompositions:
\begin{align*}
    \Sh(l_{\alpha_1+\alpha_2} l _{\alpha_1+\alpha_3}) &= (l_{\alpha_1+\alpha_2}, l_{\alpha_1+\alpha_3}),
\\    \Sh(l_{\alpha_1+\alpha_2} l _{\alpha_2+\alpha_3}) &= (l_{\alpha_1+\alpha_2}, l_{\alpha_2+\alpha_3}),
\\    \Sh(l_{\alpha_1+\alpha_3} l _{\alpha_2+\alpha_3}) &= (l_{\alpha_1+\alpha_3}, l_{\alpha_2+\alpha_3}),
\end{align*}
which give place to relations \eqref{genqSerre-ej2} by Corollary \ref{Coro:otrasrelaciones}.

The decomposition $\Sh(l_{\alpha_1} l _{\alpha_2+\alpha_3}) = (l_{\alpha_1}, l_{\alpha_2+\alpha_3})$ tell us that
$[ x_1, x_{\alpha_2+\alpha_3}]_c$ is a linear combination of $x_{\alpha_1+\alpha_2+\alpha_3}$ and $x_2 x_{\alpha_1+\alpha_3}$
by Lemma \ref{Lemma:otrasrelaciones}, and we calculate the corresponding coefficients using Lemma
\ref{Lemma:calculo coef}.

Also $\Sh(l_{\alpha_1} l _{\alpha_1+\alpha_2+\alpha_3})= (l_{\alpha_1}, l_{\alpha_1+\alpha_2+\alpha_3})$, so
$$ [x_1, x_{\alpha_1+\alpha_2+\alpha_3}]_c=0. $$
This relation is again redundant because of \eqref{idjac}, $x_1^2=0$ and the first relation in \eqref{genqSerre-ej2}.
The same holds for the relation $[x_{\alpha_1+\alpha_2+\alpha_3}, x_2]_c= 0$,
coming from the decomposition $\Sh(l _{\alpha_1+\alpha_2+\alpha_3}l_{\alpha_2}) = (l_{\alpha_1+\alpha_2+\alpha_3},
l_{\alpha_2})$.

Also $\Sh(l_{\alpha_1+\alpha_2} l _{\alpha_1+\alpha_2+\alpha_3})= (l_{\alpha_1+\alpha_2}, l_{\alpha_1+\alpha_2+\alpha_3})$, so
$$ [ x_{\alpha_1+\alpha_2}, x_{\alpha_1+\alpha_2+\alpha_3}]_c = 0. $$
This relation is also redundant by the previous relations and \eqref{idjac}. In the same way,
$ [ x_{\alpha_1+\alpha_2+\alpha_3}, x_{\alpha_1+\alpha_3}]_c = [ x_{\alpha_1+\alpha_2+\alpha_3},
x_{\alpha_2+\alpha_3} ]_c = 0 $ are redundant. The proposition follows by Theorem \ref{Thm:presentacion}.
\ep

\begin{obs}\label{Rem:generacion en rango uno}
We can prove that if $(V,c)$ is a braided vector space as in Proposition \ref{Prop:ejemplo1} or Proposition \ref{Prop:ejemplo2}, and $R= \oplus_{n \geq 0} R_n$ is a finite-dimensional graded braided Hopf algebra such that $R_0= \kk 1$ and $R_1 \cong V$ as braided vector spaces, then $R$ is generated by $R_1$ as an algebra. The proof is exactly as in \cite[Thm. 2.7]{AnGa}, using the corresponding presentation by generators and relations.
\end{obs}

\begin{obs}\label{Rem:standard}
When the braiding is of standard type, we obtain the presentation by generators and relations given in \cite[Section 5]{An}. In fact, Corollary \ref{Coro:caracterizacionLyndon} gives the set of Lyndon words obtained in \cite[Section 4B]{An}. Then we obtain a set of relations as in Theorem \ref{Thm:presentacion}, where the set of relations \eqref{quantumSerregeneralizadas1} includes the ones \cite[Theorems 5.14, 5.19, 5.22, 5.25]{An} which are not root vectors powers. Then we can reduce this set of relations because of \eqref{idjac} as in such paper, in order to obtain a minimal set of relations.
\end{obs}

\medskip

\subsection*{Acknowledgements}

I want to thank my advisor Nicol\'as Andruskiewitsch for his excellent guidance, his suggestions and his reading of
this work.

\end{document}